\documentclass[12pt,leqno,twoside]{article} 
\usepackage{bezier,amsmath,amssymb,stmaryrd,lscape,mathdots}

\input xy
\xyoption{all}
\input{rk.sty}

\renewcommand{\baselinestretch}{1.1}
\newcommand{\ulk}[1]{\ul{#1\!}\,}
\DeclareMathOperator{\per}{periodic}
\DeclareMathOperator{\TTT}{{\sf T}}
\DeclareMathOperator{\III}{{\sf I}}
\newcommand{\hsm}{\hspace*{-1mm}}
\newcommand{\tru}{\vartriangle}
\newcommand{\trd}{\triangledown}
\newcommand{\dDe}{\dot\De}
\newcommand{\ulEl}{\ul{\El\!}\,}
\newcommand{\uulEl}{\ulk{\smash{\ulEl}\rule[-0.2ex]{0mm}{0mm}}}
\newcommand{\trud}{{\tru\!\trd}}

\begin{document}

\title{Nonisomorphic Verdier octahedra on the same base}
\author{Matthias K\"unzer}
\maketitle

\begin{small}
\begin{quote}
\begin{center}{\bf Abstract}\end{center}\vspace*{2mm}
We show by an example that in a Verdier triangulated category, there may exist two mutually nonisomorphic Verdier octahedra containing the same commutative triangle.
\end{quote}
\end{small}

\renewcommand{\thefootnote}{\fnsymbol{footnote}}
\renewcommand{\thefootnote}{\arabic{footnote}}
\newdir{ >}{{}*!/-5pt/@{>}}

\begin{footnotesize}
\renewcommand{\baselinestretch}{0.7}%
\parskip0.0ex%
\tableofcontents%
\parskip1.2ex%
\renewcommand{\baselinestretch}{1.0}%
\end{footnotesize}%

\setcounter{section}{-1}

\section{Introduction}

\subsection{Is being a $3$-triangle characterised by $2$-triangles?}
\label{SecIntroDV}

{\sc Verdier} (implicitly) defined a Verdier octahedron to be a diagram in a triangulated category in the shape of an octahedron, four of whose triangles are distinguished, the 
four others commutative \bfcite{Ve63}{Def.\ 1-1}; cf.\ also \bfcite{BBD84}{1.1.6}. It arises as follows. 

To a morphism in a triangulated category, we can attach an object, called its {\it cone}. The morphism we start with and its cone are contained in a distinguished triangle. To the morphism we started 
with, we refer as the {\it base} of this distinguished triangle.

Now given a commutative triangle, we can form the cone on the first morphism, on the second morphism and on their composite, yielding three distinguished triangles. These three cones in turn are
contained in a fourth distinguished triangle. The whole diagram obtained by this construction is a Verdier octahedron. We shall refer to the commutative triangle we started with as the {\it base} of this 
Verdier octahedron.

A distinguished triangle has the property of being determined up to isomorphism by its base. Moreover, any morphism between the bases of two distinguished triangles
can be extended to a morphism between the whole distinguished triangles.

We shall show that the analogous assertion is not true for Verdier octahedra. In \S\ref{SecNonisoVO}, we give an example of two nonisomorphic Verdier octahedra on the same base. In particular, 
the identity morphism between the bases cannot be prolonged to a morphism between the whole Verdier octahedra.

The reader particularly interested in Verdier octahedra can read \S\ref{SecBaseCat}, \S\ref{SecShiftEx}, \S\ref{SecVerdierEx} and \S\ref{SecNonisoVO}.

In the terminology of Heller triangulated categories, a Verdier octahedron is a periodic $3$\nbd-pretriangle $X$ such that $Xd^\#$ is a $2$\nbd-triangle (i.e.\ a distinguished triangle) 
for all injective periodic monotone maps $\b\De_3\llaa{d} \b\De_2$.

One of the two Verdier octahedra in our example will be a $3$\nbd-triangle in the sense of \mb{\bfcite{Ku05}{Def.\ 1.5}}, i.e.\ a ``distinguished octahedron'', whereas the other will not.

Note that unlike a Verdier octahedron, a $3$\nbd-triangle is uniquely determined up to isomorphism by its base in the Heller triangulated context; cf.~\bfcite{Ku05}{Lem.\ 3.4.(6)}.

\subsection{Is being an $n$-triangle characterised by $(n-1)$-triangles?}

The situation of \S\ref{SecIntroDV} can be generalised in the following manner. 

Suppose given a closed Heller triangulated category $(\Cl,\TTT,\tht)$; cf.~\bfcite{Ku05}{Def.\ 1.5}, Definition~\ref{DefClosed}. 

The Heller triangulation $\tht = (\tht_n)_{n\ge 0}$ on $(\Cl,\TTT)$ can be viewed as a means to distinguish certain periodic $n$\nbd-pretriangles as $n$\nbd-triangles. Namely, 
a periodic $n$\nbd-pretriangle $X$ is, by definition, an $n$\nbd-triangle if $X\tht_n = 1$; cf.~\bfcite{Ku05}{Def.\ 1.5.(ii.2)}. For instance, $2$\nbd-triangles are distinguished triangles in 
the sense of Verdier; $3$\nbd-triangles are particular, ``distinguished'' Verdier octahedra.

\subsubsection{The example}

Let $n\ge 3$. Let $X$ be a periodic $n$\nbd-pretriangle. Suppose that $Xd^\#$ is an $(n-1)$\nbd-triangle for all injective periodic monotone maps $\b\De_n\llaa{d} \b\De_{n-1}$. One might ask whether 
$X$ is an $n$\nbd-triangle.

We shall show in \S\ref{SecNonisoPre} by an example that this is, in general, not the case. 

\subsubsection{Consequences}

Suppose given $n\ge 3$ and a subset of the set of periodic $n$\nbd-pretriangles. We shall say for the moment that {\it determination} holds for this subset if for $X$ and $\w X$ 
out of this subset, $X|_{\dDe_n}\iso\w X|_{\dDe_n}\ru{-1.5}$ implies that there is a periodic isomorphism $X\iso\w X$. We shall say that {\it prolongation} holds for this subset, if for $X$ and $\w X$ 
out of this subset and a morphism $X|_{\dDe_n}\lra\w X|_{\dDe_n}$, there exists a periodic morphism $X\lra\w X$ that restricts on $\dDe_n$ to that given morphism. If prolongation holds, then 
determination holds. 

\begin{itemize}
\item Consider the subset of periodic $n$\nbd-pretriangles $X$ such that $Xd^\#$ is an $(n-1)$\nbd-triangle for all injective periodic monotone maps $\b\De_n\llaa{d} \b\De_{n-1}$. Our example shows 
that in general, determination and prolongation do not hold for this subset. In fact, if $X$ is such an $n$\nbd-pretriangle, but not an $n$\nbd-triangle, then the $n$\nbd-triangle on the base $X|_{\dDe_n}$ is 
not isomorphic to $X$; cf.~\bfcite{Ku05}{Lem.~3.4.(1,\,4)}.

\item {\sc Bernstein}, {\sc Beilinson} and {\sc Deligne} considered the subset of periodic $n$\nbd-pretriangles $X$ such that $Xd^\#$ is a $2$\nbd-triangle (i.e.\ a distinguished triangle) for all 
injective periodic monotone maps $\b\De_n\llaa{d} \b\De_2$ \bfcite{BBD84}{1.1.14}. Our \mb{example} shows that in general, determination and prolongation do not hold for this subset. In fact, this 
subset contains the previously described subset.
\end{itemize}

In both of the cases above, if $n = 3$, then the condition singles out the subset of Verdier octahedra. 

\begin{itemize}
\item By \bfcite{Ku05}{Lem.\ 3.4.(6); Lem.\ 3.2}, determination and prolongation hold for the set of $n$\nbd-triangles. 
\end{itemize}

So morally, our example shows that it makes sense to let the Heller triangulation $\tht$ distinguish $n$\nbd-triangles for all $n\ge 0$. There is no ``sufficiently large'' $n$ we could be content with. 

\subsection{An appendix on transport of structure}

Suppose given a Frobenius category $\El$; that is, an exact category with enough bijective objects (relative to pure short exact sequences). Let $\Bl\tm\El$ denote the full subcategory of bijective objects.

There are two variants of the stable category of $\El$. First, there is the {\it classical stable category} $\ulEl$, defined as the quotient of $\El$ modulo $\Bl$. Second, there is the 
{\it stable category} $\uulEl$, defined as the quotient of the category of purely acyclic complexes with entries in $\Bl$ modulo the category of split acyclic complexes with entries in $\Bl$.
The categories $\ulEl$ and $\uulEl$ are equivalent. The advantage of the variant $\uulEl$ is that it carries a shift automorphism, whereas $\ulEl$ carries a shift autoequivalence. 

In \bfcite{Ku05}{Cor.\ 4.7}, we have endowed $\uulEl$ with a Heller triangulation. Now in our particular situation, also $\ulEl$ carries a shift automorphism. Since $\ulEl$ is better suited for 
calculations within that category, the question arises whether the equivalence $\ulEl \iso \uulEl$ can be used to transport the structure of a Heller triangulated category from $\uulEl$ to $\ulEl$. 
This is indeed the case; cf.\ Proposition~\ref{PropTDS}.(1). Moreover, we give recipes how to detect and how to construct $n$\nbd-triangles in $\ulEl$; cf.\ 
Propositions~\ref{PropTDS}.(2,\,3),~\ref{PropTDS7}.

Roughly put, the variant $\uulEl$ is rather suited for theoretical purposes, the variant $\ulEl$ is rather suited for practical purposes, and we had to pass a result from $\uulEl$ to $\ulEl$. Not 
surprisingly, to do so, we had to grapple with the various equivalences and isomorphisms involved.

\subsection{Acknowledgements}

I thank {\sc Amnon Neeman} for pointing out, years ago, why a counterexample as in \S\ref{SecNonisoPre} should exist, contrary to what I had believed.

This example has been found using the computer algebra system {\sc Magma}~\bfcit{CaEtAl}. I thank {\sc Markus Kirschmer} for help with a Magma program. 

I thank the referee for helpful comments.

\subsection{Notations and conventions} 
\label{SecNotConv}

\begin{footnotesize}
We use the conventions listed in \bfcite{Ku05}{\S 0.6}. In addition, we use the following conventions.

\begin{itemize}
\item[(i)] If $x$ and $y$ are elements of a set, we let $\dell_{x,y} := 1$ if $x = y$, and we let $\dell_{x,y} := 0$ if $x \ne y$.
\item[(ii)] Given $a\in\Z$, we write $\Z/a := \Z/a\Z$.
\item[(iii)] Given a ring $R$ and $R$\nbd-modules $X$ and $Y$, we write, by choice, $\liu{R}{(X,Y)} = \liu{R\Modl}{(X,Y)} = \Hom_R(X,Y)$. Moreover, given $k\ge 0$, we write 
$X^{\ds k} := \Ds_{i\in [1,k]} X$.
\item[(iv)] An {\it automorphism} $T$ of a category $\Cl$ is an endofunctor on $\Cl$ for which there exists an endofunctor $S$ such that $ST = 1_\Cl$ and $TS = 1_\Cl$. 
An {\it autoequivalence} $T$ of a category $\Cl$ is an endofunctor on $\Cl$ for which there exists an endofunctor $S$ such that $ST \iso 1_\Cl$ and $TS \iso 1_\Cl$.
\item[(v)] Let $n\ge 0$. Recall that $\b\De_n^\trud = \{\be/\al\in\b\De_n^\# : 0 \le \al\le \be\le 0^{+1}\} \tm\b\De_n^\#$.
We will often display an $n$\nbd-triangle or a periodic $n$\nbd-pretriangle in a Heller triangulated category $\Cl$ by showing its restriction to 
$\b\De_n^\trud \ohne (\{ \al/\al\; :\; 0\le\al\le 0^{+1}\} \cup \{ 0^{+1}/0 \})$. This is possible without loss of information, for we can reconstruct the whole diagram by adding zeroes on $\al/\al$
for $0\le\al\le 0^{+1}$ and on $0^{+1}/0$, and then by periodic prolongation.
\item[(vi)] Suppose given a Heller triangulated category $\Cl$. A {\it Verdier octahedron} in $\Cl$ is a periodic $3$\nbd-pretriangle $X\in\Ob\,\Cl^{+,\,\per}(\b\De_3^\#)$ such that 
$X d^\#\in\Ob\,\Cl^{+,\,\per}(\b\De_2^\#)$ is a $2$\nbd-triangle for all injective periodic monotone maps $\smash{\b\De_3\llaa{d} \b\De_2}\ru{4.5}$.
\end{itemize}
\end{footnotesize}

\fbox{Henceforth, let $p\ge 2$ be a prime.}

\section{The classical stable category of $(\Z/p^m)\modl$}
\label{SecExZmodpm}

\subsection{The category $(\Z/p^m)\modl$}
\label{SecBaseCat}

Let $m\ge 0$. By $\El := (\Z/p^m)\modl$ we understand the following category. 

The objects are indexed by tuples $(a_i)_{i\in [0,m]}$ with $a_i\in\Z_{\ge 0}$. To such an index, we attach the object
\[
\Ds_{i\in [0,m]} (\Z/p^i)^{\ds a_i}\; .
\]
As morphisms, we take $\Z/p^m$\nbd-linear maps. 

\bq
 Note that we have {\sf not} chosen a skeleton. The trick here is to pick {\sf several} zero objects.
\eq

The duality contrafunctor $\liu{\Z/p^m}{(-,\Z/p^m)}$ on $\El$, which sends $\Z/p^i$ to $\Z/p^i$ for $i\in [1,m]$, shows that an object in this category is injective 
if and only if it is projective. An object of $\El$ is bijective if and only if it is isomorphic to a finite direct sum of copies of $\Z/p^m$.
The category $\El$ is an abelian Frobenius category, with all short exact sequences stipulated to be pure; cf.\ e.g.\ \mb{\bfcite{Ku05}{Def.\ A.5.(2)}}.

\subsection{The shift on \underline{$(\Z/p^m)\modl$}}
\label{SecShiftEx}

To define a shift automorphism on the classical stable category $\ulEl = \ul{(\Z/p^m)\modl}$, we shall distinguish certain (pure) short exact sequences in $\El$; 
cf.~\S\ref{SecHypo},~\bfcite{Ku05}{Def.\ A.7}.

Let $\EE_k := \smatdd{1}{}{\vspace*{-2mm}}{}{\ddots}{}{}{}{1}$ denote the unit matrix of size $k\ti k$; let $\EE'_k := \smatdd{}{}{1\vspace*{-2mm}}{}{\iddots}{}{1}{}{}$ denote the reversed unit 
matrix of size $k\ti k$.

As distinguished (pure) short exact sequences we take those of the form
\[
\hspace*{-5mm}
\xymatrix@C=40mm{
\Ds_{i\in [0,m]} (\Z/p^i)^{\ds a_i} \arm[r]^{\enger{\left(\ba{cccc} \scm p^m\EE_{a_0}\vspace*{2mm} &&& \\ &\scm p^{m-1}\EE_{a_1} && \\ && \ddots & \\ &&&\scm p^0\EE_{a_m} \\ \ea\right)}} & 
(\Z/p^m)^{\ds\sum_{i\in [0,m]} a_i} \are[r]^{\enger{\left(\ba{cccc}  &&&\scm\EE'_{a_0}\vspace*{2mm} \\ &&\scm\EE'_{a_1}\vspace*{-2mm} & \\ & \iddots && \\ \scm \EE'_{a_m} &&& \\ \ea\right)}} &
\Ds_{i\in [0,m]} (\Z/p^i)^{\ds a_{m-i}}
}
\]
So roughly speaking, distinguished short exact sequences are direct sums of those of the form
\[
\Z/p^i\;\;\;\lramonoa{p^{m-i}}\;\;\; \Z/p^m\;\;\;\lraepia{1}\;\;\; \Z/p^{m-i}\; ,
\]
where $i\in [0,m]$; we reorder the summands the cokernel term consists of.

With this choice, conditions (i,\,ii,\,iii) of \S\ref{SecHypo} are satisfied.

On indecomposable objects and morphisms between them, the shift automorphism induced on $\ulEl$ by our set of distinguished short exact sequences is given  by 
\[
(\Z/p^i\;\lrafl{25}{a}\;\Z/p^j)^{+1} \;\=\;  (\Z/p^{m-i}\;\lrafl{25}{p^{i-j} a}\;\Z/p^{m-j})\;,
\]
where $i,\, j\,\in\,[0,m]$, and where $a$ is a representative in $\Z$. Note that if $i < j$, then $a$ is divisible by $p^{\,j-i}$.

Note that $\Z/p^i\lraa{a}\Z/p^j$ represents zero in $\ulEl$ if and only if $a$ is divisible by $p^{\min(m-i,\, j)}$.

\subsection{A Heller triangulation on \underline{$(\Z/p^m)\modl$}}
\label{SecHellerEx}

Concerning the notation $\El^\Box(\b\De_n^\trud)$, cf.~\S\ref{SecSomeLem}. Given $n\ge 0$ and $X\in\Ob\El^\Box(\b\De_n^\trud)$, we form $X^\tau\in\Ob\ulEl^{+,\,\per}(\b\De_n^\trud)$ with 
respect to the set of distinguished short exact sequences of \S\ref{SecShiftEx} as described in \S\ref{SecStandardise}. That is, we replace the rightmost column of $X$ by the 
column obtained using distinguished short exact sequences, so that $(X^\tau)_{0^{+1}/\ast} = ((X^\tau)_{\ast/0})^{+1} = (X_{\ast/0})^{+1}$; cf.~\S\ref{SecStandardise}.

\begin{Remark}
\label{RemDist}
If the short exact sequences 
\[
X_{\al/0}\;\lramonofl{25}{\smatez{x}{x}}\; X_{\al/\al}\ds X_{0^{+1}/0}\;\lraepifl{40}{\rsmatze{x}{-x}}\; X_{0^{+1}/\al}
\]
appearing in the diagram $X$ for $1\le\al\le n$ already are distinguished, then the image of $X$ in $\Ob\ulEl^+(\bar\De_n^\trud)$ equals $X^\tau$.
\end{Remark}

Concerning the notion of a closed Heller triangulated category, cf.\ Definition~\ref{DefClosed} in \S\ref{SecDetecting}.

\begin{Remark}
\label{RemEx}
The classical stable category $\,\ulEl = \ul{(\Z/p^m)\modl}\ru{-2.6}\,$ carries a closed Heller triangulation such that given $n\ge 0$ and $X\in\Ob\El^\Box(\b\De_n^\trud)$, the periodic prolongation of 
$X^\tau$ to an object of $\ulEl^{+,\,\per}(\b\De_n^\#)$ is an $n$\nbd-triangle.
\end{Remark}

{\it Proof.} The assertion follows by Proposition~\ref{PropTDS}.(1) in \S\ref{SecHypo} and Proposition~\ref{PropTDS7} in \S\ref{SecTrianglesClassical}. \qed

\subsection{A Verdier triangulation on \underline{$(\Z/p^m)\modl$}} 
\label{SecVerdierEx}

By \bfcite{Ha88}{Th.\ 2.6}, $\ulEl = \ul{(\Z/p^m)\modl}$ is a Verdier triangulated category, i.e.\ a triangulated category in the sense of {\sc Verdier} \bfcite{Ve63}{Def.\ 1-1}. 

\bq
 This also follows by Remark~\ref{RemEx} and by \bfcite{Ku05}{Prop.\ 3.6}, which says that any Heller triangulated category in which idempotents split is also Verdier triangulated. The $2$\nbd-triangles 
 in the Heller context are the distinguished triangles in the Verdier context.
\eq

Given a morphism $X\lraa{f} Y$ in $\El$, using the distinguished short exact sequence $X\lramono B\lraepi X^{+1}$, where $B$ is bijective, we can form the morphism
\[
\xymatrix{
Y\arm[r]         & Z\arefl{0.4}[r]       & X^{+1}           \\
X\arm[r]\ar[u]^f & B\arefl{0.4}[r]\ar[u] & X^{+1}\ar@{=}[u] \\
}
\]
of short exact sequences, from which the sequence
\[
\xymatrix{
X\ar[r]^f & Y\arm[r] & Z\arefl{0.4}[r] & X^{+1} \\
}
\]
represents a distinguished triangle in the Verdier triangulated category $\ulEl$.

\section{Nonisomorphic periodic $n$-pretriangles}
\label{SecNonisoPre}

\bq
 Nonisomorphic periodic $n$\nbd-pretriangles whose periodic $(n-1)$\nbd-pretriangles are all $(n-1)$\nbd-triangles, to be specific.
\eq

Let $n\ge 3$. Let $\Cl := \ul{(\Z/p^{2n})\modl}$, and let it be endowed with a shift automorphism as in \S\ref{SecShiftEx} and a Heller triangulation as in \S\ref{SecHellerEx}.

\subsection{A $(2n-1)$-triangle}
\label{SecInflatedTriangle}

Let $Y$ be the following $(2n-1)$\nbd-triangle in $\Cl$. 

{\footnotesize
\[
\xymatrix{
               &                        &                                             &                     &                             &                                & \Z/p^1                \\
               &                        &                                             &                     &                             & \Z/p^1\ar[r]^{-p}              & \Z/p^2\ar[u]^1        \\
               &                        &                                             &                     & \Z/p^1\ar[r]^p              & \Z/p^2\ar[r]^{-p}\ar[u]^1      & \Z/p^3\ar[u]^1        \\
               &                        & \hspace*{16mm}\iddots\ru{-2}\hspace*{-16mm} &                     & \vdots\ru{-2}\ar[u]^1       & \vdots\ru{-2}\ar[u]^1          & \vdots\ru{-2}\ar[u]^1 \\
               &                        & \Z/p^1\ar[r]^p                              & \cdots\ar[r]^(0.4)p & \Z/p^{2n-4}\ar[r]^p\ar[u]^1 & \Z/p^{2n-3}\ar[r]^{-p}\ar[u]^1 & \Z/p^{2n-2}\ar[u]^1   \\
               & \Z/p^1\ar[r]^p         & \Z/p^2\ar[r]^p\ar[u]^1                      & \cdots\ar[r]^(0.4)p & \Z/p^{2n-3}\ar[r]^p\ar[u]^1 & \Z/p^{2n-2}\ar[r]^{-p}\ar[u]^1 & \Z/p^{2n-1}\ar[u]^1   \\
\Z/p^1\ar[r]^p & \Z/p^2\ar[r]^p\ar[u]^1 & \Z/p^3\ar[r]^p\ar[u]^1                      & \cdots\ar[r]^(0.4)p & \Z/p^{2n-2}\ar[r]^p\ar[u]^1 & \Z/p^{2n-1}\ar[u]^1            &                       \\
}
\]
}

\bq
 Here we have made use of the convention from \S\ref{SecNotConv} that we display of $Y$ only its restriction to the subposet 
 $\{\be/\al\in\b\De_{2n-1}^\#\; :\; 0\le\al < \be \le 0^{+1},\; \be/\al \ne 0^{+1}/0\}$, which is possible without loss of information. Similarly below.
\eq

It arises from a diagram on $\bar\De_n^\trud$ with values in $(\Z/p^{2n})\modl$ that consists of squares, has entry $\Z/p^{2n}$ at position $0^{+1}/0$, and has the quadrangle
\[
\xymatrix{
\Z/p^{2n-2}\ar[r]^{-p}      & \Z/p^{2n-1}          \\
\Z/p^{2n-1}\ar[r]^p\ar[u]^1 & \Z/p^{2n}\ar[u]_{-1} \\
}
\]
in its lower right corner. This diagram contains the necessary distinguished short exact sequences with the necessary signs inserted for $Y$ to be in fact a 
$(2n-1)$\nbd-triangle; cf.\ Remarks~\ref{RemDist},~\ref{RemEx}.

\subsection{An $n$-triangle and a periodic $n$-pretriangle}
\label{SecTriangle}

We apply the folding operator $\ffk_{n-1}$ to the $(2n-1)$\nbd-triangle $Y$ obtained in \S\ref{SecInflatedTriangle}, yielding the $n$\nbd-triangle $Y\ffk_{n-1}$, which we shall display now; 
cf.~\bfcite{Ku05}{Lem.\ 3.4.(2), \S 1.2.2.3}.

{\footnotesize
\[
\xymatrix@C=6mm{
                 &                                                  &                                                                                   &                                         &                                                                                 &                                                                                                   & \Z/p^n               \\
                 &                                                  &                                                                                   &                                         &                                                                                 & \Z/p^1\dk\Z/p^{2n-1}\ar[r]^(0.65){\rsmatze{\hsm-p^{n-1}\hsm}{-1\;}}                               & \Z/p^n\ar[u]_p       \\
                 &                                                  &                                                                                   &                                         & \Z/p^1\dk\Z/p^{2n-1}\ar[r]^{\smatzz{p}{0}{0}{1}}                                & \Z/p^2\dk\Z/p^{2n-2}\ar[r]^(0.65){\rsmatze{\hsm-p^{n-2}\hsm}{-1\;}}\ar[u]^{\smatzz{1}{0}{0}{p}}   & \Z/p^n\ar[u]_p       \\
                 &                                                  & \hspace*{18mm}\iddots\ru{-2}\hspace*{-18mm}                                       &                                         & \vdots\ru{-2}\ar[u]^{\smatzz{1}{0}{0}{p}}                                       & \vdots\ru{-2}\ar[u]^{\smatzz{1}{0}{0}{p}}                                                         & \vdots\ru{-2}\ar[u]_p\\
                 &                                                  & \Z/p^1\dk\Z/p^{2n-1}\ar[r]^(0.7){\smatzz{p}{0}{0}{1}}                             & \cdots\ar[r]^(0.3){\smatzz{p}{0}{0}{1}} & \Z/p^{n-3}\dk\Z/p^{n+3}\ar[r]^{\smatzz{p}{0}{0}{1}}\ar[u]^{\smatzz{1}{0}{0}{p}} & \Z/p^{n-2}\dk\Z/p^{n+2}\ar[r]^(0.65){\rsmatze{\hsm -p^2\hsm}{-1}}\ar[u]^{\smatzz{1}{0}{0}{p}}     & \Z/p^n\ar[u]_p       \\
                 & \Z/p^1\dk\Z/p^{2n-1}\ar[r]^{\smatzz{p}{0}{0}{1}} & \Z/p^2\dk\Z/p^{2n-2}\ar[r]^(0.7){\smatzz{p}{0}{0}{1}}\ar[u]^{\smatzz{1}{0}{0}{p}} & \cdots\ar[r]^(0.3){\smatzz{p}{0}{0}{1}} & \Z/p^{n-2}\dk\Z/p^{n+2}\ar[r]^{\smatzz{p}{0}{0}{1}}\ar[u]^{\smatzz{1}{0}{0}{p}} & \Z/p^{n-1}\dk\Z/p^{n+1}\ar[r]^(0.65){\rsmatze{-p}{-1}}\ar[u]^{\smatzz{1}{0}{0}{p}}                & \Z/p^n\ar[u]_p       \\
\Z/p^n\ar[r]^{p} & \Z/p^n\ar[r]^{p}\ar[u]_{\smatez{1}{-p^{n-1}}}    & \Z/p^n\ar[r]^{p}\ar[u]_{\smatez{1}{-p^{n-2}}}                                     & \cdots\ar[r]^{p}                        & \Z/p^n\ar[r]^{p}\ar[u]_{\smatez{1}{-p^2}}                                       & \Z/p^n\ar[u]_{\smatez{1}{-p}}                                                                     &                      \\
} 
\]
}

Let $X$ be the $n$\nbd-triangle obtained from $Y\ffk_{n-1}$ by isomorphic substitution along $\rsmatzz{1}{0}{0}{-1}$ on all terms consisting of two summands; cf. \bfcite{Ku05}{Lem.\ 3.4.(4)}. So $X$ can
be displayed as follows.

{\footnotesize
\[
\xymatrix@C=6mm{
                 &                                                  &                                                                                   &                                         &                                                                                 &                                                                                                    & \Z/p^n               \\
                 &                                                  &                                                                                   &                                         &                                                                                 & \Z/p^1\dk\Z/p^{2n-1}\ar[r]^(0.65){\rsmatze{\hsm-p^{n-1}\hsm}{1\;\;\;}}                             & \Z/p^n\ar[u]_p       \\
                 &                                                  &                                                                                   &                                         & \Z/p^1\dk\Z/p^{2n-1}\ar[r]^{\smatzz{p}{0}{0}{1}}                                & \Z/p^2\dk\Z/p^{2n-2}\ar[r]^(0.65){\rsmatze{\hsm-p^{n-2}\hsm}{1\;\;\;}}\ar[u]^{\smatzz{1}{0}{0}{p}} & \Z/p^n\ar[u]_p       \\
                 &                                                  & \hspace*{18mm}\iddots\ru{-2}\hspace*{-18mm}                                       &                                         & \vdots\ru{-2}\ar[u]^{\smatzz{1}{0}{0}{p}}                                       & \vdots\ru{-2}\ar[u]^{\smatzz{1}{0}{0}{p}}                                                          & \vdots\ru{-2}\ar[u]_p\\
                 &                                                  & \Z/p^1\dk\Z/p^{2n-1}\ar[r]^(0.7){\smatzz{p}{0}{0}{1}}                             & \cdots\ar[r]^(0.3){\smatzz{p}{0}{0}{1}} & \Z/p^{n-3}\dk\Z/p^{n+3}\ar[r]^{\smatzz{p}{0}{0}{1}}\ar[u]^{\smatzz{1}{0}{0}{p}} & \Z/p^{n-2}\dk\Z/p^{n+2}\ar[r]^(0.65){\rsmatze{\hsm -p^2\hsm}{1\;}}\ar[u]^{\smatzz{1}{0}{0}{p}}     & \Z/p^n\ar[u]_p       \\
                 & \Z/p^1\dk\Z/p^{2n-1}\ar[r]^{\smatzz{p}{0}{0}{1}} & \Z/p^2\dk\Z/p^{2n-2}\ar[r]^(0.7){\smatzz{p}{0}{0}{1}}\ar[u]^{\smatzz{1}{0}{0}{p}} & \cdots\ar[r]^(0.3){\smatzz{p}{0}{0}{1}} & \Z/p^{n-2}\dk\Z/p^{n+2}\ar[r]^{\smatzz{p}{0}{0}{1}}\ar[u]^{\smatzz{1}{0}{0}{p}} & \Z/p^{n-1}\dk\Z/p^{n+1}\ar[r]^(0.65){\rsmatze{-p}{1}}\ar[u]^{\smatzz{1}{0}{0}{p}}                  & \Z/p^n\ar[u]_p       \\
\Z/p^n\ar[r]^{p} & \Z/p^n\ar[r]^{p}\ar[u]_{\smatez{1}{p^{n-1}}}     & \Z/p^n\ar[r]^{p}\ar[u]_{\smatez{1}{p^{n-2}}}                                      & \cdots\ar[r]^{p}                        & \Z/p^n\ar[r]^{p}\ar[u]_{\smatez{1}{p^2}}                                        & \Z/p^n\ar[u]_{\smatez{1}{p}}                                                                       &                      \\
} 
\]
}

Let $\w X$ be the following periodic $n$\nbd-pretriangle.

{\footnotesize
\[
\xymatrix@C=6mm{
                  &                                                  &                                                                                   &                                         &                                                                                             &                                                                                                     & \Z/p^n               \\
                  &                                                  &                                                                                   &                                         &                                                                                             & \Z/p^1\dk\Z/p^{2n-1}\ar[r]^(0.65){\rsmatze{\hsm -p^{n-1}\hsm}{1\;\;\;}}                             & \Z/p^n\ar[u]_p       \\
                  &                                                  &                                                                                   &                                         & \Z/p^1\dk\Z/p^{2n-1}\ar[r]^{\smatzz{p}{0}{0}{1}}                                            & \Z/p^2\dk\Z/p^{2n-2}\ar[r]^(0.65){\rsmatze{\hsm -p^{n-2}\hsm}{1\;\;\;}}\ar[u]^{\smatzz{1}{0}{0}{p}} & \Z/p^n\ar[u]_p       \\
                  &                                                  & \hspace*{18mm}\iddots\ru{-2}\hspace*{-18mm}                                       &                                         & \vdots\ru{-2}\ar[u]^{\smatzz{1}{0}{0}{p}}                                                   & \vdots\ru{-2}\ar[u]^{\smatzz{1}{0}{0}{p}}                                                           & \vdots\ru{-2}\ar[u]_p\\
                  &                                                  & \Z/p^1\dk\Z/p^{2n-1}\ar[r]^(0.7){\smatzz{p}{0}{0}{1}}                             & \cdots\ar[r]^(0.3){\smatzz{p}{0}{0}{1}} & \Z/p^{n-3}\dk\Z/p^{n+3}\ar[r]^{\smatzz{p}{0}{0}{1}}\ar[u]^{\smatzz{1}{0}{0}{p}}             & \Z/p^{n-2}\dk\Z/p^{n+2}\ar[r]^(0.65){\rsmatze{\hsm -p^2\hsm}{1\;}}\ar[u]^{\smatzz{1}{0}{0}{p}}      & \Z/p^n\ar[u]_p       \\
                  & \Z/p^1\dk\Z/p^{2n-1}\ar[r]^{\smatzz{p}{0}{0}{1}} & \Z/p^2\dk\Z/p^{2n-2}\ar[r]^(0.7){\smatzz{p}{0}{0}{1}}\ar[u]^{\smatzz{1}{0}{0}{p}} & \cdots\ar[r]^(0.3){\smatzz{p}{0}{0}{1}} & \Z/p^{n-2}\dk\Z/p^{n+2}\ar[r]_(0.45){\smatzz{p}{0}{p^{n-3}}{1}}\ar[u]^{\smatzz{1}{0}{0}{p}} & \Z/p^{n-1}\dk\Z/p^{n+1}\ar[r]^(0.65){\rsmatze{-p}{1}}\ar[u]^{\smatzz{1}{0}{-p^{n-3}}{p}}            & \Z/p^n\ar[u]_p       \\
\Z/p^n\ar[r]^{p}  & \Z/p^n\ar[r]^{p}\ar[u]_{\smatez{1}{p^{n-1}}}     & \Z/p^n\ar[r]^{p}\ar[u]_{\smatez{1}{p^{n-2}}}                                      & \cdots\ar[r]^{p}                        & \Z/p^n\ar[r]^{p}\ar[u]_{\smatez{1}{p^2}}                                                    & \Z/p^n\ar[u]_{\smatez{1}{p}}                                                                        &                      \\
} 
\]
}

To verify that $\w X$ actually is an $n$\nbd-pretriangle, a comparison with $X$ reduces us to show that the three quadrangles depicted in full in the lower right corner of $\w X$ are weak squares.
Of these three, the middle quadrangle arises from the corresponding one of $X$ by an isomorphic substitution along 
$\Z/p^{n-1}\ds\Z/p^{n+1}\;\;\lraisofl{45}{\smatzz{1}{0}{p^{n-3}}{1}}\;\;\Z/p^{n-1}\ds\Z/p^{n+1}\ru{9}$, and thus is a weak square. For the lower one, we may apply \bfcite{Ku05}{Lem.\ A.17}
to the diagram $(\w X_{1/0},\,\w X_{n-1/0},\,\w X_{n/0},\,\w X_{0^{+1}/0},\,\w X_{1/1},\,\w X_{n-1/1},\,\w X_{n/1},\,\w X_{0^{+1}/1})$ and compare with $X$ to show that it is a weak square.
For the right hand side one, we may apply \bfcite{Ku05}{Lem.\ A.17} to the diagram 
$(\w X_{n/0},\,\w X_{n/1},\,\w X_{n/2},\,\w X_{n/n},\,\w X_{0^{+1}/0},\,\w X_{0^{+1}/1},\,\w X_{0^{+1}/2},\,\w X_{0^{+1}/n})$ and compare with $X$ to show that it is a weak square.

Given $k\in [0,n]$, we let $\b\De_n\llafl{25}{\dd_k}\b\De_{n-1}$ be the periodic monotone map determined by \mb{$[0,n-1]\dd_k = [0,n]\ohne\{ k\}$}.

\begin{Lemma}
\label{LemIsCovered} 
Suppose given $k\in [0,n]$.
\begin{itemize}
\item[{\rm (1)}] The diagram $\w X \dd_k^\#$ is an $(n-1)$\nbd-triangle. 
\item[{\rm (2)}] We have $X \dd_k^\# \iso \w X \dd_k^\#$ in $\Cl^{+,\,\per}(\b\De_{n-1}^\#)$. 
\end{itemize}
\end{Lemma}

{\it Proof.} Since $X$ is an $n$\nbd-triangle, $X \dd_k^\#$ is an $(n-1)$\nbd-triangle; cf.~\bfcite{Ku05}{Lem.\ 3.4.(1)}. Since $X \dd_k^\#|_{\dDe_{n-1}} = \w X \dd_k^\#|_{\dDe_{n-1}}$, the diagram 
$\w X \dd_k^\#$ is an $(n-1)$\nbd-triangle if and only if it is isomorphic to $X \dd_k^\#$ in $\Cl^{+,\,\per}(\b\De_{n-1}^\#)$; cf.~\bfcite{Ku05}{Lem.\ 3.4.(4,\,6)}. So assertions (1) and (2) are 
equivalent. We will prove (2).

When referring to an object on a certain position in the diagram $X \dd_k^\#$ resp.\ $\w X \dd_k^\#$, we shall also mention in parentheses its position as an object in the diagram $X$ resp.\ $\w X$ for 
ease of orientation.

When constructing a morphism in $\Cl^{+,\,\per}(\b\De_{n-1}^\#)$, we will give its components on \mb{$\{ j/i \; :\; 0\le i\le j\le n-1\}\tm \b\De_{n-1}^\#\,$}; the remaining components result 
thereof by periodic repetition.

{\it Case $k\in\{ 1,\, n\}$}. We have $X \dd_1^\# = \w X \dd_1^\#$ and $X \dd_n^\# = \w X \dd_n^\#$.

{\it Case $k = 0$.} We {\it claim} that $X \dd_0^\#$ is isomorphic to $\w X \dd_0^\#$ in $\Cl^{+,\,\per}(\b\De_{n-1}^\#)$. In fact, an isomorphism $X \dd_0^\#\lraiso\w X \dd_0^\#$ is given by 
\[
\Z/p^{n-1}\ds\Z/p^{n+1}\;\;\mraisofl{43}{\smatzz{1}{0}{p^{n-3}}{1}}\;\;\Z/p^{n-1}\ds\Z/p^{n+1}
\]
at position $(n-1)/0$ (position $n/1$ in $X$ resp.\ $\w X$), and by the identity elsewhere. This proves the {\it claim.}

{\it Case $k\in [2,n-1]$.} We {\it claim} that $X \dd_k^\#$ is isomorphic to $\w X \dd_k^\#$ in $\Cl^{+,\,\per}(\b\De_{n-1}^\#)$. In fact, an isomorphism $X \dd_k^\#\lraiso\w X \dd_k^\#$ is given as follows.

At position $j/0$ for $j\in [1,n-1]$ (position $j/0$ if $j\le k-1$ and $(j+1)/0$ if $j\ge k$ in $X$ resp.\ $\w X$), it is given by the identity on $\Z/p^n$. 

At position $j/i$ for $i,\, j\,\in\, [1,k-1]$ such that $i < j$ (position $j/i$ in $X$ resp.\ $\w X$), it is given by the identity on
$\Z/p^{\,j-i}\ds\Z/p^{2n-j+i}$.

At position $j/i$ for $i,\, j\,\in\, [k,n-1]$ such that $i < j$ (position $(j+1)/(i+1)$ in $X$ resp.\ $\w X$), it is given by the identity on
$\Z/p^{\,j-i}\ds\Z/p^{2n-j+i}$.

At position $j/i$ for $i\in [1,k-1]$ and $j\in [k,n-1]$ such that $j/i\ne (n-1)/1$ (position $(j+1)/i$ in $X$ resp.\ $\w X$), it is 
given by
\[
\Z/p^{\,j+1-i}\ds\Z/p^{2n-j-1+i}\;\;\mraisofl{45}{\smatzz{1}{0}{-p^{\,j-1-i}}{1}}\;\;\Z/p^{\,j+1-i}\ds\Z/p^{2n-j-1+i}\ru{7}
\]

At position $(n-1)/1$ (position $n/1$ in $X$ resp.\ $\w X$), it is given by the identity on $\Z/p^{n-1}\ds\Z/p^{n+1}$. 

This proves the {\it claim.}\qed

\begin{Lemma}
\label{LemNonIso}
$X$ is not isomorphic to $\w X$ in $\Cl^{+,\,\per}(\b\De_n^\#)$.
\end{Lemma}

In particular, $\w X$ is not an $n$\nbd-triangle; cf.~\bfcite{Ku05}{Lem.\ 3.4.(6)}.

{\it Proof.} We {\it assume} the contrary. By \bfcite{Ku05}{3.4.(4)}, $X$ and $\w X$ are $n$\nbd-triangles. Thus, by \bfcite{Ku05}{3.4.(6)}, there is an isomorphism $X\lraiso\w X$ that is identical 
at $i/0$ and at $0^{+1}/i$ for $i\in [1,n]$. Let 
\[
\Z/p^{\ell - k}\dk\Z/p^{2n - \ell + k} \;\;\vlrafl{43}{\smatzz{\;a_{\ell/k}\;}{p^{2n-2\ell+2k} b_{\ell/k}}{c_{\ell/k}}{d_{\ell/k}}}\ru{7}\;\; \Z/p^{\ell - k}\dk\Z/p^{2n - \ell + k}
\]
denote the entry of this isomorphism at $\ell/k$, where $1\le k < \ell \le n$.

If $\ell - k \ge 2$, we have the following commutative quadrangle in $\Cl$ on $\ell/k \lra \ell/(k+1)$.
\[
\xymatrix@C=40mm@R=15mm{
\Z/p^{\ell - k}\dk\Z/p^{2n - \ell + k}\ar[d]_{\smatzz{\;a_{\ell/k}\;}{p^{2n-2\ell+2k} b_{\ell/k}}{c_{\ell/k}}{d_{\ell/k}}} \ar[r]^{\smatzz{1}{0}{0}{p}} 
                                                  & \Z/p^{\ell - k - 1}\dk\Z/p^{2n - \ell + k + 1}\ar[d]^{\smatzz{\;a_{\ell/(k+1)}\;}{p^{2n-2\ell+2k+2} b_{\ell/(k+1)}}{c_{\ell/(k+1)}}{d_{\ell/(k+1)}}} \\
\Z/p^{\ell - k}\dk\Z/p^{2n - \ell + k}\ar[r]_(0.47){\smatzz{1}{0}{0}{p} \;-\; \dell_{\ell/k,\,n/1}\smatzz{0}{0}{p^{n-3}}{0}} 
                                                  & \Z/p^{\ell - k - 1}\dk\Z/p^{2n - \ell + k + 1} \\ 
}
\]
We read off the congruences
\begin{align}
c_{\ell/k} - \dell_{\ell/k,\,n/1}\, p^{n-3} d_{\ell/k} & \;\;\con_{p^{\ell-k-1}}\;\;  p c_{\ell/(k+1)} \tag{i}  \\
b_{\ell/k}                                             & \;\;\con_{p^{\ell-k-1}}\;\;  p b_{\ell/(k+1)} \; . \tag{ii} 
\end{align}

From (i) we infer
\[
c_{n/1} - p^{n-3} d_{n/1} \;\con_{p^{n-2}}\; p^1 c_{n/2} \;\con_{p^{n-2}}\; p^2 c_{n/3} \;\con_{p^{n-2}}\; \dots \;\con_{p^{n-2}}\; p^{n-2} c_{n/(n-1)} \;\con_{p^{n-2}}\; 0\; . \tag{iii}
\]
From (ii) we infer
\[
b_{n/1} \;\con_{p^{n-2}}\; p b_{n/2} \;\con_{p^{n-2}}\; p^2 b_{n/3} \;\con_{p^{n-2}}\; \dots \;\con_{p^{n-2}}\; p^{n-2} b_{n/(n-1)} \;\con_{p^{n-2}}\; 0\; . \tag{iv}
\]

On $n/1\lra 0^{+1}/1$, we have the following commutative quadrangle in $\Cl$.
\[
\xymatrix{
\Z/p^{n-1}\dk\Z/p^{n+1}\ar[d]_{\smatzz{a_{n/1}}{p^2 b_{n/1}}{c_{n/1}}{d_{n/1}}}\ar[r]^(0.65){\rsmatze{-p}{1}} & \Z/p^n\ar@{=}[d] \\
\Z/p^{n-1}\dk\Z/p^{n+1}\ar[r]_(0.65){\rsmatze{-p}{1}}                                                         & \Z/p^n           \\
}
\]
We read off the congruence
\[
- p c_{n/1} + d_{n/1} \;\;\con_{p^{n-1}}\;\; 1 \; . \tag{v}
\]
On $n/0\lra n/1$, we have the following commutative quadrangle in $\Cl$.
\[
\xymatrix{
\Z/p^n\ar@{=}[d]\ar[r]^(0.3){\smatez{1}{p}} & \Z/p^{n-1}\dk\Z/p^{n+1}\ar[d]^{\smatzz{a_{n/1}}{p^2 b_{n/1}}{c_{n/1}}{d_{n/1}}} \\
\Z/p^n          \ar[r]_(0.3){\smatez{1}{p}} & \Z/p^{n-1}\dk\Z/p^{n+1}                                                         \\
}
\]
We read off the congruence
\[
p b_{n/1} + d_{n/1} \;\;\con_{p^{n-1}}\;\; 1 \; . \tag{vi}
\]

By (iii) resp.\ (iv) we conclude from (v) resp.\ (vi) that
\begin{align}
(1 - p^{n-2})d_{n/1} &\;\;\con_{p^{n-1}}\;\; 1      \tag{v$'$} \\
d_{n/1}              &\;\;\con_{p^{n-1}}\;\; 1 \; . \tag{vi$'$}   
\end{align}
Substituting (vi$'$) into (v$'$), we obtain
\[
1 - p^{n-2}\;\;\con_{p^{n-1}}\;\; 1\; ,
\]
which is {\it absurd.}\qed

\section{Nonisomorphic Verdier octahedra}
\label{SecNonisoVO}

\bq
 Since in \S\ref{SecNonisoPre}, the category $\Cl$ is also a Verdier triangulated category, specialising to $n = 3$ yields two nonisomorphic Verdier octahedra on the same base. In this particular 
 case, we shall now give a somewhat longer argument alternative to that given in \S\ref{SecNonisoPre} that is independent of \bfcit{Ku05}, whose techniques 
 might not be familiar to all readers. Nonetheless, \S\ref{SecNonisoVO} is a particular case of \S\ref{SecNonisoPre}.
\eq

Let $\Cl := \ul{(\Z/p^6)\modl}$, and let it be endowed with a shift automorphism as in \S\ref{SecShiftEx} and a Verdier triangulation as in \S\ref{SecVerdierEx}.

Let the diagram $X$ be given by 
\[
\xymatrix{
               &                                             &                                                                           & \Z/p^3         \\
               &                                             & \Z/p^1\dk\Z/p^5\ar[r]^(0.60){\rsmatze{-p^2}{1\;}}                         & \Z/p^3\ar[u]_p \\
               & \Z/p^1\dk\Z/p^5\ar[r]^{\smatzz{p}{0}{0}{1}} & \Z/p^2\dk\Z/p^4\ar[r]^(0.60){\rsmatze{-p}{1}}\ar[u]^{\smatzz{1}{0}{0}{p}} & \Z/p^3\ar[u]_p \\
\Z/p^3\ar[r]^p & \Z/p^3\ar[r]^p\ar[u]_{\smatez{1}{p^2}}      & \Z/p^3\ar[u]_{\smatez{1}{p}}                                              & .              \\
}
\]

Let the diagram $\w X$ be given by
\[
\xymatrix{
               &                                             &                                                                             & \Z/p^3         \\
               &                                             & \Z/p^1\dk\Z/p^5\ar[r]^(0.60){\rsmatze{-p^2}{1\;}}                           & \Z/p^3\ar[u]_p \\
               & \Z/p^1\dk\Z/p^5\ar[r]^{\smatzz{p}{0}{1}{1}} & \Z/p^2\dk\Z/p^4\ar[r]^(0.60){\rsmatze{-p}{1}}\ar[u]^{\rsmatzz{1}{0}{-1}{p}} & \Z/p^3\ar[u]_p \\
\Z/p^3\ar[r]^p & \Z/p^3\ar[r]^p\ar[u]_{\smatez{1}{p^2}}      & \Z/p^3\ar[u]_{\smatez{1}{p}}                                                & .              \\
}
\]

\begin{Lemma}
\label{LemVerdierOct}
Both $X$ and $\w X$ are Verdier octahedra. 
\end{Lemma}

\bq
 In contrast to the procedure in \S\ref{SecNonisoPre}, to prove this, we will not make use of the folding operation. 
\eq

{\it Proof.} For the periodic monotone map $\b\De_3\llafl{25}{\dd_3}\b\De_2$ that maps $0\lamaps 0$, $1\lamaps 1$ and $2\lamaps 2$, we obtain $X \dd_3^\# = \w X \dd_3^\#$, horizontally displayed as
\[
\xymatrix{
\Z/p^3\ar[r]^p & \Z/p^3\ar[r]^(0.4){\smatez{1}{p^2}} & \Z/p\dk\Z/p^5\ar[r]^(0.6){\rsmatze{-p^2}{1\;}} & \Z/p^3 \; . \\
}
\] 
The following morphism of short exact sequences in $(\Z/p^6)\modl$ shows $X \dd_3^\#$ to be a distinguished triangle.
\[
\xymatrix{
\Z/p^3\armfl{0.4}[r]^(0.4){\smatez{1}{p^2}} & \Z/p\dk\Z/p^5\arefl{0.6}[r]^(0.6){\rsmatze{-p^2}{1\;}}              & \Z/p^3           \\
\Z/p^3\armfl{0.4}[r]^(0.4){p^3}\ar[u]^p     & \Z/p^6       \arefl{0.6}[r]^(0.6){1\ru{-0.5}}\ar[u]_{\smatez{0}{1}} & \Z/p^3\ar@{=}[u] \\
}
\]

For the periodic monotone map $\b\De_3\llafl{25}{\dd_1}\b\De_2$ that maps $0\lamaps 0$, $2\lamaps 1$ and $3\lamaps 2$, we obtain the distinguished triangle $X \dd_1^\# = \w X \dd_1^\# = X \dd_3^\#$
again. 

For the periodic monotone map $\b\De_3\llafl{25}{\dd_2}\b\De_2$ that maps $0\lamaps 0$, $1\lamaps 1$ and $3\lamaps 2$, we obtain the diagram $X \dd_2^\# = \w X \dd_2^\#$, horizontally displayed as
\[
\xymatrix{
\Z/p^3\ar[r]^{p^2} & \Z/p^3\ar[r]^(0.4){\smatez{1}{p}} & \Z/p^2\dk\Z/p^4\ar[r]^(0.6){\rsmatze{-p}{1\;}} & \Z/p^3 \; . \\
}
\] 
The following morphism of short exact sequences in $(\Z/p^6)\modl$ shows $X \dd_2^\#$ to be a distinguished triangle.
\[
\xymatrix{
\Z/p^3\armfl{0.4}[r]^(0.4){\smatez{1}{p}}   & \Z/p^2\dk\Z/p^4\arefl{0.6}[r]^(0.6){\rsmatze{-p}{1\;}}                & \Z/p^3           \\
\Z/p^3\armfl{0.4}[r]^(0.4){p^3}\ar[u]^{p^2} & \Z/p^6         \arefl{0.6}[r]^(0.6){1\ru{-0.5}}\ar[u]_{\smatez{0}{1}} & \Z/p^3\ar@{=}[u] \\
}
\]
For the periodic monotone map $\b\De_3\llafl{25}{\dd_0}\b\De_2$ that maps $1\lamaps 0$, $2\lamaps 1$ and $3\lamaps 2$, we obtain the periodic isomorphism $X \dd_0^\# \lraiso \w X \dd_0^\#$, horizontally 
displayed as 
\[
\xymatrix@C=12mm{
\Z/p\dk\Z/p^5\ar[r]^{\smatzz{p}{0}{0}{1}}\ar@{=}[d] & \Z/p^2\dk\Z/p^4\ar[r]^{\smatzz{1}{0}{0}{p}}\ar[d]^{\smatzz{1}{0}{1}{1}}_\wr & \Z/p\dk\Z/p^5\ar[r]^{\smatzz{0}{-p^4}{1}{0}}\ar@{=}[d] & \Z/p\dk\Z/p^5\ar@{=}[d] \\
\Z/p\dk\Z/p^5\ar[r]^{\smatzz{p}{0}{1}{1}}           & \Z/p^2\dk\Z/p^4\ar[r]^{\rsmatzz{1}{0}{-1}{p}}                               & \Z/p\dk\Z/p^5\ar[r]^{\smatzz{0}{-p^4}{1}{0}}           & \Z/p\dk\Z/p^5\zw{.}     \\
}
\] 
So we are reduced to show that $X \dd_0^\#$ is a distinguished triangle, which it is as a direct sum of two distinguished triangles, as the following morphisms of short exact sequences in $(\Z/p^6)\modl$ show.
\[
\xymatrix{
\Z/p  \arm[r]^p         & \Z/p^2\are[r]^{1\ru{-0.5}}         & \Z/p           & & \Z/p^4\arm[r]^p         & \Z/p^5\are[r]^{1\ru{-0.5}}         & \Z/p           \\
\Z/p^5\arm[r]^p\ar[u]^1 & \Z/p^6\are[r]^{1\ru{-0.5}}\ar[u]_1 & \Z/p\ar@{=}[u] & & \Z/p^5\arm[r]^p\ar[u]^1 & \Z/p^6\are[r]^{1\ru{-0.5}}\ar[u]_1 & \Z/p\ar@{=}[u] \\
}
\]
\qed

\begin{Lemma}
\label{LemOctNonIso}
The Verdier octahedra $X$ and $\w X$ are not isomorphic in $\Cl^{+,\,\per}(\b\De_3^\#)$.
\end{Lemma}

That is, there is no isomorphism between the displayed parts of $X$ resp.\ of $\w X$ such that its entries on the rightmost vertical column arise by an application of the shift functor of $\Cl$ to its
entries on the lower row.

\bq
 We will not use the fact that $X$ is a $3$\nbd-triangle, which in conjunction with \bfcite{Ku05}{3.4.(4,\,6)} would permit us to restrict ourselves to consider isomorphisms that are identical 
 on the lower row and the rightmost vertical column, as we did in Lemma~\ref{LemNonIso}.
\eq

{\it Proof.} We {\it assume} the contrary and depict an isomorphism $X\lraiso\w X$ as follows.
{\footnotesize
\[
\hspace*{-5mm}
\xymatrix@!@C=-2mm@R=-5mm{
                          & &                                                                     &                                                                                            &                                                          &                                                                                                                           &                                                                                                &                                &                           & \Z/p^3\ar[ddd]^w \\ 
                          & &                                                                     &                                                                                            &                                                          &                                                                                                                           & \Z/p\dk\Z/p^5\ar[rr]^{\rsmatze{-p^2}{1\;}}\ar'[d][ddd]^(0.35){\smatzz{a''}{p^4 b''}{c''}{d''}} &                                & \Z/p^3\ar[ur]^p\ar[ddd]^v &                  \\
                          & &                                                                     & \Z/p\dk\Z/p^5\ar[rr]^{\smatzz{p}{0}{0}{1}}\ar'[d][ddd]^(0.35){\smatzz{a'}{p^4 b'}{c'}{d'}} &                                                          & \Z/p^2\dk\Z/p^4\ar[ur]^*+<-2mm,0mm>{\smatzz{1}{0}{0}{p}}\ar[rr]^(0.68){\rsmatze{-p}{1}}\ar[ddd]^{\smatzz{a}{p^2 b}{c}{d}} &                                                                                                & \Z/p^3\ar[ur]^p\ar[ddd]^(0.4)u &                           &                  \\
\Z/p^3\ar[rr]^p\ar[ddd]^u & & \Z/p^3\ar[rr]^(0.6)p\ar[ur]^*+<-2mm,0mm>{\smatez{1}{p^2}}\ar[ddd]^v &                                                                                            & \Z/p^3\ar[ur]^*+<-2mm,0mm>{\smatez{1}{p}}\ar[ddd]^(0.4)w &                                                                                                                           &                                                                                                &                                &                           & \Z/p^3           \\  
                          & &                                                                     &                                                                                            &                                                          &                                                                                                                           & \Z/p\dk\Z/p^5\ar'[r][rr]^(0.4){\rsmatze{-p^2}{1\;}}                                            &                                & \Z/p^3\ar[ur]^p           &                  \\
                          & &                                                                     & \Z/p\dk\Z/p^5\ar'[r]^(0.72){\smatzz{p}{0}{1}{1}}[rr]                                       &                                                          & \Z/p^2\dk\Z/p^4\ar[ur]^*+<-4mm,0mm>{\rsmatzz{1}{0}{-1}{p}}\ar[rr]^{\rsmatze{-p}{1}}                                       &                                                                                                & \Z/p^3\ar[ur]^p                &                           &                  \\
\Z/p^3\ar[rr]^p           & & \Z/p^3\ar[rr]^p\ar[ur]^*+<-4mm,0mm>{\smatez{1}{p^2}}                &                                                                                            & \Z/p^3\ar[ur]^*+<-4mm,0mm>{\smatez{1}{p}}                &                                                                                                                           &                                                                                                &                                &                           &                  \\  
}
\]
}

Note that all vertical quadrangles commute in $\Cl$.

The commutative quadrangles on $1/0\lra 2/0\lra 3/0$ yield $u\con_{p^2} v\con_{p^2} w$.

The commutative quadrangle on $3/0\lra 3/1$ yields $pb + d\con_{p^2} w$.

The commutative quadrangle on $3/1\lra 0^{+1}/1$ yields $-pc + d\con_{p^2} u$.

The commutative quadrangle on $3/1\lra 3/2$ yields $b\con_p 0$ and $c\con_p d$.

Altogether, we have 
\[
u\;\con_{p^2}\; w \;\con_{p^2}\; pb + d \;\con_{p^2}\; d \;\con_{p^2}\; u + pc \;\con_{p^2}\; u + pd\; ,
\]
whence 
\[
0\;\con_p\; d\;\con_p\; w\; .
\]
Since $\Z/p^3\lraa{w}\Z/p^3$ is an isomorphism in $\Cl$, we have $w\not\con_p 0$. This is {\it absurd.} \qed

In \bfcite{BBD84}{1.1.13}, it is described how an octahedron gives rise to two ``extra'' triangles. As cone of the diagonal of a quadrangle appearing in that octahedron, we take the direct sum 
of the non-diagonal terms of the subsequent quadrangle, the morphisms being taken from the octahedron, with one minus sign inserted to ensure that the composition of two morphisms in the constructed 
triangle vanishes.

\begin{Remark}
\label{RemBBD}
The triangles arising from $X$ and from $\w X$ as described in \bfcite{BBD84}{1.1.13} are distinguished.
\end{Remark}

{\it Proof.} (\footnote{Strictly speaking, we should reorder summands in the diagrams that follow; cf.~\S\ref{SecBaseCat}. But then the proof would be more difficult to read.}) 
The morphism of short exact sequences in $(\Z/p^6)\modl$ 
\[
\xymatrix{
\Z/p^3\armfl{0.34}[r]^(0.34){\smated{1}{\;p^2\;}{-p}}                                        & \Z/p\ds\Z/p^5\ds\Z/p^3\arefl{0.55}[r]^(0.55){\smatdz{p}{0}{0}{1}{1}{p}}                           & \Z/p^2\ds\Z/p^4           \\
\Z/p^4\ds\Z/p^2\armfl{0.47}[r]^{\smatzz{p^2}{0}{0}{p^4}}\ar[u]_{\smatze{-p\ru{-1.3}}{\;p^2}} & \Z/p^6\ds\Z/p^6\arefl{0.50}[r]^(0.50){\smatzz{1}{0}{0}{1}}\ar[u]_{\rsmatzd{0}{-p}{\; 1}{0}{1}{0}} & \Z/p^2\ds\Z/p^4\ar@{=}[u] \\
}
\]
and the isomorphism of diagrams with coefficients in $\Cl$
\[
\xymatrix@C=20mm{
\Z/p^3\ar[r]^(0.35){\smated{1\;}{p^2}{-p}}\ar@{=}[d] & \Z/p\dk\Z/p^5\dk\Z/p^3\ar[r]^(0.55){\smatdz{p}{0}{0}{1}{1}{p}}\ar@{=}[d] & \Z/p^2\dk\Z/p^4\ar[r]^(0.55){\rsmatze{-p^2}{p\;}}\ar[d]^{\smatzz{1-p\;}{0}{1\;}{1}}_\wr & \Z/p^3\ar@{=}[d] \\
\Z/p^3\ar[r]^(0.35){\smated{1\;}{p^2}{-p}}           & (\Z/p\dk\Z/p^5)\dk\Z/p^3\ar[r]^(0.55){\smatdz{p}{0}{1}{1}{1}{p}}         & \Z/p^2\dk\Z/p^4\ar[r]^(0.55){\rsmatze{-p^2}{p\;}}                                       & \Z/p^3           \\
}
\]
show one of the triangles mentioned in loc.\ cit.\ to be distinguished in $X$ and in $\w X$.

The morphism of short exact sequences in $(\Z/p^6)\modl$
\[
\xymatrix@C=12mm{
\Z/p^2\ds\Z/p^4\armfl{0.42}[r]^(0.42){\rsmatzd{1}{\;0}{p}{0}{p}{-1}} & \Z/p\ds\Z/p^5\ds\Z/p^3\arefl{0.67}[r]^(0.67){\rsmatde{-p^2}{1\;}{p\;}} & \Z/p^3           \\
\Z/p^3\arm[r]^{p^3}\ar[u]_{\smatez{p}{p^2}}                          & \Z/p^6\are[r]^{1\ru{-1}}\ar[u]_{\smated{0}{1}{0}}                      & \Z/p^3\ar@{=}[u] \\
}
\]
and the isomorphism of diagrams with coefficients in $\Cl$
\[
\xymatrix@C=23mm@R=12mm{
\Z/p^3\ar[r]^(0.4){\smatez{p\;}{p^2}}\ar@{=}[d] & \Z/p^2\dk\Z/p^4\ar[r]^(0.43){\rsmatzd{1}{\;0}{p}{0}{p}{-1}}\ar@{=}[d] & (\Z/p\dk\Z/p^5)\dk\Z/p^3\ar[r]^(0.6){\rsmatde{-p^2}{1\;}{p\;}}\ar[d]^{\smatdd{1}{\;0}{-p^2}{0}{1}{1}{1}{0}{1+p}}_\wr & \Z/p^3\ar@{=}[d] \\
\Z/p^3\ar[r]^(0.4){\smatez{p\;}{p^2}}           & \Z/p^2\dk\Z/p^4\ar[r]^(0.43){\rsmatzd{1}{\;0}{p}{-1}{p}{-1}}          & (\Z/p\dk\Z/p^5)\dk\Z/p^3\ar[r]^(0.6){\rsmatde{-p^2}{1\;}{p\;}}                                                       & \Z/p^3           \\
}
\]
show the other of the triangles mentioned in loc.\ cit.\ to be distinguished in $X$ and in $\w X$.
\qed

\appendix

\begin{footnotesize}

\section{Transport of structure}

We use the notation of \bfcite{Ku05}{\S1,\,\S2}.

\subsection{Transport of a Heller triangulation}

Concerning weakly abelian categories, see e.g.\ \bfcite{Ku05}{\S A.6.3}.
Recall that an additive functor between weakly abelian categories is called subexact if it induces an exact functor on the Freyd categories; cf.~\bfcite{Ku05}{\S 1.2.1.3}. For instance, an equivalence is 
subexact.

\begin{Setup}
\label{Setup1}\rm\Absit

\begin{tabular}{p{1mm}p{15cm}}
& Suppose given a Heller triangulated category $(\Cl,\TTT,\tht)$; cf.~\bfcite{Ku05}{Def.\ 1.5}. Suppose given a weakly abelian category $\Cl'$ and an automorphism $\TTT'$ on 
$\Cl'$, called {\it shift}\,; cf.~\bfcite{Ku05}{Def.\ A.26}. \\
&  Assume given subexact functors $\Cl\lraa{F}\Cl'$ and $\Cl'\lraa{G}\Cl$, and isotransformations $1_{\Cl'}\lraisoa{\eps} GF$ and $\smash{\TTT' G\lraisoa{\sa} G\TTT}\ru{4}$. \\
\end{tabular}
\end{Setup}

Suppose given $n\ge 0$. By abuse of notation, we write $F := \ulk{F^+(\b\De_n^\#)} : \ulk{\Cl^+(\b\De_n^\#)} \lra \ulk{\Cl'^+(\b\De_n^\#)}\ru{-2.2}$ for the functor obtained by pointwise 
application of $F$. 

Similarly, we write $\eps := \ulk{\eps^+(\b\De_n^\#)} : \ulk{(1_{\Cl'})^+(\b\De_n^\#)} \lraiso \ulk{(GF)^+(\b\De_n^\#)}$ for the isotransformation obtained by pointwise application of $\eps$. 

More generally speaking, for notational convenience, induced functors of type $\ulk{A^+(\b\De_n^\#)}$ will often be abbreviated by $A$, and induced transformations of type $\ulk{\al^+(\b\De_n^\#)}$ will 
often be abbreviated by $\al$. For instance, given \mb{$X\in\Ob\,\ulk{\Cl^+(\b\De_n^\#)}$}, we will allow ourselves to write $X\TTT = X\ulk{\TTT^+(\b\De_n^\#)} \;(\; = [X^{+1}])$.

Given $X'\in\Ob\,\ulk{\Cl'^+(\b\De_n^\#)} = \Ob\,\Cl'^+(\b\De_n^\#)$, we define the isomorphism $[X']^{+1}\lraisofl{25}{X'\tht'_n} [X'^{+1}]$ in $\ulk{\Cl'^+(\b\De_n^\#)}$ by the following commutative 
diagram.
\[
\xymatrix@C=16mm{
[X']^{+1}\ar[r]^{X'\tht'_n}_\sim\ar[d]_{[X'\eps]^{+1}}^\wr & [X'^{+1}]\ar[d]^{[X'^{+1}]\eps}_\wr \\
[X'GF]^{+1}\ar@{=}[d]                                      & [X'^{+1}]GF\ar[d]^\wr_{X'\sa F}     \\
[X'G]^{+1} F\ar[r]^{X'G\tht_n F}_\sim                      & [(X' G)^{+1}]F                      \\
}
\]
In other words, we let 
\[
X'\tht'_n \; :=\; ([X'\eps]^{+1})(X'G\tht_n F)(X'\sa^- F)([X'^{+1}]\eps^-)\; .
\]
As a composite of isotransformations, $(X'\tht'_n)_{X'\,\in\,\Ob\,\ulk{\Cl'^+(\b\De_n^\#)}}$ is an isotransformation. Let $\tht' := (\tht'_n)_{n\ge 0}$.

\begin{Lemma}
\label{LemTDS1}
The triple $(\Cl',\TTT',\tht')$ is a Heller triangulated category.
\end{Lemma}

Cf.~\bfcite{Ku05}{Def.\ 1.5}. We will say that $\tht'$ is {\it transported from $(\Cl,\TTT,\tht)$ via $F$ and $G$.} Strictly speaking, we should mention $\eps$ and $\sa$ here as well.

{\it Proof.} Suppose given $m,\, n\,\ge\, 0$, a periodic monotone map $\b\De_n\llaa{q}\b\De_m$ and $X'\in\Ob\,\ulk{\Cl'^+(\b\De_n^\#)}$. We {\it claim} that 
$X'\ulk{q}^\#\tht'_m = X'\tht'_n\ulk{q}^\#$. We have
\[
\ba{lcl}
X'\ulk{q}^\#\tht'_m & = & ([X'\ulk{q}^\#\eps]^{+1}) (X'\ulk{q}^\# G\tht_m F) (X'\ulk{q}^\#\sa^- F) ([{X'\ulk{q}^\#}^{+1}]\eps^-)\; ,\\
X'\tht'_n\ulk{q}^\# & = & ([X'\eps]^{+1}\ulk{q}^\#) (X'G\tht_n F\ulk{q}^\#)  (X'\sa^- F\ulk{q}^\#) ([X'^{+1}]\eps^-\ulk{q}^\#)\; . \\
\ea
\]
By respective pointwise definition, we have $[X'\ulk{q}^\#\eps]^{+1} = [X'\eps]^{+1}\ulk{q}^\#$ (using that $q$ is periodic), $X'\ulk{q}^\#\sa^- F = X'\sa^- F\ulk{q}^\#$ and 
$[{X'\ulk{q}^\#}^{+1}]\eps^- = [X'^{+1}]\eps^-\ulk{q}^\#$. Moreover, since $(\Cl,\TTT,\tht)$ is Heller triangulated, we get
\[
X'\ulk{q}^\# G\tht_m F\= X'G\ulk{q}^\#\tht_m F \= X'G\tht_n\ulk{q}^\# F \= X'G\tht_n F\ulk{q}^\#\; .
\]
This proves the {\it claim.}

Suppose given $n\ge 0$ and $X'\in\Ob\,\ulk{\Cl'^+(\b\De_{2n+1}^\#)}$. We {\it claim} that $X'\ul{\ffk}_n\tht'_{n+1} = X'\tht'_{2n+1}\ul{\ffk}_n$. We have
\[
\ba{lcl}
X'\ul{\ffk}_n\tht'_{n+1}  & = & ([X'\ul{\ffk}_n\eps]^{+1}) (X'\ul{\ffk}_n G\tht_{n+1} F) (X'\ul{\ffk}_n\sa^- F) ([{X'\ul{\ffk}_n}^{+1}]\eps^-)\; ,\\
X'\tht'_{2n+1}\ul{\ffk}_n & = & ([X'\eps]^{+1}\ul{\ffk}_n) (X'G\tht_{2n+1} F\ul{\ffk}_n)  (X'\sa^- F\ul{\ffk}_n) ([X'^{+1}]\eps^-\ul{\ffk}_n)\; . \\
\ea
\]
By additivity of $F$, $G$ and $\TTT'$ and by respective pointwise definition, we have $[X'\ul{\ffk}_n\eps]^{+1} = [X'\eps]^{+1}\ul{\ffk}_n$ (using shiftcompatibility of $\ul{\ffk}_n$), 
$X'\ul{\ffk}_n\sa^- F = X'\sa^- F\ul{\ffk}_n$ and $[{X'\ul{\ffk}_n}^{+1}]\eps^- = [X'^{+1}]\eps^-\ul{\ffk}_n$. Moreover, since $(\Cl,\TTT,\tht)$ is Heller triangulated, we get
\[
X'\ul{\ffk}_n G\tht_{n+1} F\= X'G\ul{\ffk}_n\tht_{n+1} F \= X'G\tht_{2n+1}\ul{\ffk}_n F \= X'G\tht_{2n+1} F\ul{\ffk}_n\; .
\]
This proves the {\it claim.}\qed

\subsection{Detecting $n$-triangles}
\label{SecDetecting}

\begin{Setup}
\label{Setup2}\rm\Absit

\begin{tabular}{p{1mm}p{15cm}}
& Suppose given a Heller triangulated category $(\Cl,\TTT,\tht)$; cf.~\bfcite{Ku05}{Def.\ 1.5}. Suppose given an additive category $\Cl'$ and an automorphism $\TTT'$ on $\Cl'$, called {\it shift}. \\
& Suppose given mutually inverse equivalences $\Cl\lraa{F}\Cl'$ and $\Cl'\lraa{G}\Cl$. Note that $G\adj F$, whence there exist isotransformations $1_{\Cl'}\lraisoa{\eps} GF$ and 
  $FG\lraisoa{\et} 1_{\Cl}$ such that both $(F\eps)(\et F) = 1_F$ and $(\eps G)(G\et) = 1_G$ hold. We fix such $\eps$ and $\et$. \\
& Suppose given an isotransformation $\TTT' G\lraisoa{\sa} G\TTT$. \\
\end{tabular}
\end{Setup}

Note that $\Cl'$ is weakly abelian, being equivalent to the weakly abelian category $\Cl$.

Let $\tht'$ be transported from $(\Cl,\TTT,\tht)$ via $F$ and $G$.

That is, we let $X'\tht'_n \; :=\; ([X'\eps]^{+1})(X'G\tht_n F)(X'\sa^- F)([X'^{+1}]\eps^-)$ for $n\ge 0$ and $X'\in\Ob\,\ulk{\Cl'^+(\b\De_n^\#)}$, defining \mb{$\tht' := (\tht'_n)_{n\ge 0}$}. 

By Lemma~\ref{LemTDS1}, the triple $(\Cl',\TTT',\tht')$ is a Heller triangulated category.

Moreover, let 
\[
(\TTT F\lraisoa{\rh} F\TTT') \= (\TTT F\lraisofl{25}{\et^-\TTT F} FG\TTT F \lraisofl{25}{F\sa^- F} F\TTT' GF\lraisofl{25}{F\TTT'\eps^-} F\TTT')\; .
\]

\begin{Notation}
\label{NotTDS1_5}\rm
Suppose given $n\ge 0$. Concerning the full subposet 
\[
\b\De_n^\trud \= \{ \be/\al\in\De_n^\# : 0\le\al\le\be\le 0^{+1}\} \;\tm\; \b\De_n^\#\; , 
\]
cf.~\bfcite{Ku05}{\S2.5.1}.

\begin{itemize}
\item[(1)] Suppose given $X'\in\Ob\,\Cl'^{+,\,\per}(\b\De_n^\#)$, where periodic means $[X']^{+1} = [X'^{+1}]$; cf.~\bfcite{Ku05}{\S2.5.3}. Consider the diagram $X'G|_{\b\De_n^\trud}$. 
Denote by $X'G|_{\b\De_n^\trud}^\sa\in\Ob\,\Cl'^+(\b\De_n^\#)$ the diagram $X'G|_{\b\De_n^\trud}$ with $(X'G)_{0^{+1}/i} = X'_{0^{+1}/i} G = X'_{i/0} \TTT' G$ isomorphically replaced via 
$X'_{i/0}\,\sa$ by $X'_{i/0} G\TTT$. Denote by $X'G^\sa\in\Ob\,\Cl^{+,\,\per}(\b\De_n^\#)$ its periodic prolongation, characterised by 
$X'G^\sa|_{\b\De_n^\trud} = X'G|_{\b\De_n^\trud}^\sa$; cf.\ \bfcite{Ku05}{\S 2.5.3}. Using, for $k\ge 0$, 
\[
\left.\left([X'G]^{+k}\lraiso [X'G^{+1}]^{+(k-1)}\lraiso\cdots\lraiso [X'G^{+k}]\right)\right|_{\b\De_n^\trud}\; ,
\]
given by 
\[
X'_{(j/i)^{+k}} G \; =\; X'_{j/i}\TTT'^k G \TTT^0\;\;\lraisofl{28}{\TTT'^{k-1}\sa\TTT^0}\;\; X'_{j/i}\TTT'^{k-1} G \TTT^1\;\;\lraisofl{28}{\TTT'^{k-2}\sa\TTT^1}\;\;\cdots
\;\;\lraisofl{28}{\TTT'^0\sa\TTT^{k-1}}\;\; X'_{j/i}\TTT'^0 G \TTT^k
\]
at $j/i$ for $0\le i\le j\le 0^{+1}$, and similarly for $k\le 0$, using $\TTT'^- G\;\;\lraisofl{28}{\;\TTT'^-\sa^-\TTT^-}\; G\TTT^-$, we obtain an isomorphism $X'G\lraisoa{\phi} X'G^\sa$ in 
$\Cl^+(\b\De_n^\#)$ such that $\phi_{i/0} = 1_{X'_{i/0}G}$ and $\phi_{0^{+1}/i} = X'_{i/0}\,\sa$ for $i\in [1,n]$.

\item[(2)] Suppose given $X\in\Ob\,\Cl^{+,\,\per}(\b\De_n^\#)$. Consider the diagram $XF|_{\b\De_n^\trud}$. Denote by $XF|_{\b\De_n^\trud}^\rh\in\Ob\,\Cl^+(\b\De_n^\#)$ the diagram 
$XF|_{\b\De_n^\trud}$ with $(XF)_{0^{+1}/i} = X_{0^{+1}/i} F = X_{i/0} \TTT F$ isomorphically replaced via $X_{i/0}\,\rh$ by $X_{i/0} F\TTT'$. Denote by 
$XF^\rh\in\Ob\,\Cl'^{+,\,\per}(\b\De_n^\#)$ its periodic prolongation, characterised by $XF^\rh|_{\b\De_n^\trud} = XF|_{\b\De_n^\trud}^\rh$; cf.~\bfcite{Ku05}{\S2.5.3}.

Similarly as in (1), we have an isomorphism $XF\lraisoa{\psi} XF^\rh$ in $\Cl'^+(\b\De_n^\#)$ such that $\psi_{i/0} = 1_{X_{i/0}F}$ and $\psi_{0^{+1}/i} = X_{i/0}\,\rh$ for $i\in [1,n]$.
\end{itemize}
\end{Notation}

\begin{Lemma}
\label{LemTDS2}
Suppose given $n\ge 0$.
\begin{itemize}
\item[{\rm (1)}] Suppose given $X'\in\Ob\,\Cl'^{+,\,\per}(\b\De_n^\#)$. Then $X'$ is an $n$\nbd-triangle if and only if $X'G^\sa$ is an $n$\nbd-triangle.
\item[{\rm (2)}] Suppose given $X\in\Ob\,\Cl^{+,\,\per}(\b\De_n^\#)$. Then $X$ is an $n$\nbd-triangle if and only if $XF^\rh$ is an $n$\nbd-triangle.
\end{itemize}
\end{Lemma}

Cf.~\mb{\bfcite{Ku05}{Def.\ 1.5.(ii.2)}.}

{\it Proof.} Ad (1). Since $\tht_n$ is a transformation, there exists a commutative quadrangle
\[
\xymatrix@C=15mm{
[X'G^\sa]^{+1}\ar[r]^{X'G^\sa\tht_n}_\sim                 & [(X'G^\sa)^{+1}]                     \\
[X'G]^{+1}\ar[r]^{X'G\tht_n}_\sim\ar[u]_{[\phi]^{+1}}^\wr & [(X'G)^{+1}]\ar[u]_{[\phi^{+1}]}^\wr \\
}
\]
in $\ulk{\Cl^+(\b\De_n^\#)}\ru{-2}$. Therefore, $X'G^\sa$ is an $n$\nbd-triangle if and only if $[\phi]^{+1} = (X'G\tht_n)[\phi^{+1}]$. By \bfcite{Ku05}{Prop.\ 2.6}, this equation is equivalent to
$[\phi]^{+1}|_{\dDe_n} = (X'G\tht_n)|_{\dDe_n}[\phi^{+1}]|_{\dDe_n}$; cf.~\bfcite{Ku05}{\S2.1.1}; in other words, to
\[
X'\sa|_{\dDe_n} \= X'G\tht_n|_{\dDe_n} 
\]
as morphisms from $X'\TTT' G|_{\dDe_n}$ to $X'G\TTT|_{\dDe_n}$ in $\ulk{\Cl(\dDe_n)}\ru{-2}$. This, in turn, is equivalent to $X'\sa \= X'G\tht_n$ as morphisms from 
$X'\TTT' G$ to $X'G\TTT$ in $\ulk{\Cl^+(\b\De_n^\#)}$ by \mb{\bfcite{Ku05}{Prop.\ 2.6}}.

Now $X'$ being an $n$\nbd-triangle is equivalent to $X'\tht'_n = 1$; i.e.\ to 
\[
([X'\eps]^{+1})(X'G\tht_n F)(X'\sa^- F)([X'^{+1}]\eps^-) \= 1 \; . 
\]
Since $[X']^{+1} = [X'^{+1}]$, we have $[X'\eps]^{+1} = [X']^{+1}\eps = [X'^{+1}]\eps$, whence this equation is equivalent to $(X'G\tht_n F)(X'\sa^- F) = 1$. Since $F$ is an equivalence, this amounts to 
$(X'G\tht_n)(X'\sa^-) = 1$, as was to be shown.

Ad (2). Since $\tht'_n$ is a transformation, there exists a commutative quadrangle
\[
\xymatrix@C=15mm{
[XF^\rh]^{+1}\ar[r]^{XF^\rh\tht'_n}_\sim                 & [(XF^\rh)^{+1}]                     \\
[XF]^{+1}\ar[r]^{XF\tht'_n}_\sim\ar[u]_{[\psi]^{+1}}^\wr & [(XF)^{+1}]\ar[u]_{[\psi^{+1}]}^\wr \\
}
\]
in $\ulk{\Cl'^+(\b\De_n^\#)}\ru{-2}$. Therefore, $XF^\rh$ is an $n$\nbd-triangle if and only if $[\psi]^{+1} = (XF\tht'_n)[\psi^{+1}]$. By \bfcite{Ku05}{Prop.\ 2.6}, this equation is equivalent to
$[\psi]^{+1}|_{\dDe_n} = (XF\tht'_n)|_{\dDe_n}[\psi^{+1}]|_{\dDe_n}$; in other words, to
\[
X\rh|_{\dDe_n} \= XF\tht'_n|_{\dDe_n} 
\]
as morphisms from $X\TTT F|_{\dDe_n}$ to $X F\TTT'|_{\dDe_n}$ in $\ulk{\Cl'(\dDe_n)}\ru{-2}$. This, in turn, is equivalent to $X\rh \= XF\tht'_n$ as morphisms from 
$X\TTT F$ to $X F\TTT'$ in $\ulk{\Cl'^+(\b\De_n^\#)}$ by \mb{\bfcite{Ku05}{Prop.\ 2.6}}. Which amounts to 
\[
(X\et^-\TTT F)(XF\sa^- F)(XF\TTT'\eps^-) \= ([XF\eps]^{+1})(XFG\tht_n F)(XF\sa^- F)([XF^{+1}]\eps^-)\; ;
\]
i.e.\ to 
\[
X\et^-\TTT F \= ([X\et^- F]^{+1})(XFG\tht_n F)\; .
\]
Since $[X\et^- F]^{+1} = [X\et^-]^{+1} F$ and since $\tht_n$ is a transformation, the right hand side equals $(X\tht_n F)([(X\et^-)^{+1}]F)$, and therefore we can continue the string of equivalent 
assertions with
\[
X\et^-\TTT F \= (X\tht_n F)([(X\et^-)^{+1}]F)\; ;
\]
i.e.\ with $X\tht_n F = 1$; i.e.\ with $X\tht_n = 1$; i.e.\ with $X$ being an $n$\nbd-triangle.\qed

\begin{Definition}
\label{DefClosed}\rm
A Heller triangulated category $(\Cl,\TTT,\tht)$ is said to be {\it closed} if every morphism $X\lraa{f} Y$ therein can be completed to a $2$\nbd-triangle; i.e.\ if for all morphisms $X\lraa{f} Y$
in $\Cl$, there exists $U\in\Ob\,\Cl^{+,\,\tht=1}(\b\De_2^\#)$ with $(X\lraa{f} Y) = (U_{1/0}\lraa{u} U_{2/0})$. If this is the case, then also the Heller triangulation $\tht$ is called {\it closed.}
\end{Definition}

For instance, a Heller triangulated category whose idempotents split is closed; cf.~\bfcite{Ku05}{Prop.\ 3.6}.

Recall that $\tht'$ is transported from $(\Cl,\TTT,\tht)$ via $F$ and $G$.

\begin{Lemma}
\label{LemTDS3}
If $(\Cl,\TTT,\tht)$ is a closed Heller triangulated category, then $(\Cl',\TTT',\tht')$ is a closed Heller triangulated category.
\end{Lemma}

{\it Proof.} By Lemma~\ref{LemTDS1}, it remains to prove closedness of $(\Cl',\TTT',\tht')$. Suppose given $X'\lraa{u'} Y'$ in $\Cl'$. We have to prove that it can be prolonged to a $2$\nbd-triangle. 
Using closedness of $(\Cl,\TTT,\tht)$, we find a $2$\nbd-triangle 
\[
X'G\;\lraa{u'G}\; Y'G\;\lraa{v}\; Z\;\lraa{w}\; X'G\TTT \; .
\]
We {\it claim} that
\[
M' \;\; :=\;\; (X'\;\lraa{u'}\; Y'\;\mrafl{25}{(Y'\eps)(vF)}\; ZF\;\vlrafl{25}{(wF)(X'\sa^- F)(X'\TTT'\eps^-)}\;  X'\TTT')
\]
is a $2$\nbd-triangle in $\Cl'$. By Lemma~\ref{LemTDS2}.(1), it suffices to show that $M'G^\sa$ is a $2$\nbd-triangle in $\Cl$. Consider the periodic isomorphism with upper row $M'G^\sa$ and lower
row a $2$-triangle 
\[
\xymatrix@C=43mm{
X'G\ar[r]^{u'G}           & Y'G\ar[r]^{(Y'\eps G)(vFG)} & ZFG\ar[r]^{(wFG)(X'\sa^-FG)(X'\TTT'\eps^-G)(X'\sa)} & X'G\TTT                        \\ 
X'G\ar[r]^{u'G}\ar@{=}[u] & Y'G\ar@{=}[u]\ar[r]^v       & Z\ar[r]^w\ar[u]^\wr_{Z\et^-}                        & X'G\TTT\ar@{=}[u] \; .\!\!\!\! \\
}
\]
In fact, we have $(Y'\eps G)(vFG) = (Y'G\et^-)(vFG) = v(Z\et^-)$ and 
\[
\barcl
(Z\et^-)(wFG)(X'\sa^-FG)(X'\TTT'\eps^-G)(X'\sa) 
& = & w(X'G\TTT\et^-)(X'\sa^-FG)(X'\TTT'\eps^-G)(X'\sa) \\
& = & w (X'\sa^-)(X'\TTT'G\et^-)(X'\TTT'\eps^-G)(X'\sa) \\
& = & w (X'\sa^-)(X'\sa) \\
& = & w \; . \\
\ea
\]
This shows that $M' G^\sa$ is a $2$-triangle; cf.~\bfcite{Ku05}{Lem.\ 3.4.(4)}. This proves the {\it claim.}
\qed

\begin{Remark}
\label{RemTDS4}

Suppose given $n\ge 0$.
\begin{itemize}
\item[{\rm (1)}] Given $X\in\Ob\,\Cl^+(\b\De_n^\#)$, we have $(X\tht_n F)(X\rh) = XF\tht'_n$ in $\ulk{\Cl'^+(\b\De_n^\#)}$.
\item[{\rm (2)}] Given $X'\in\Ob\,\Cl'^+(\b\De_n^\#)$, we have $(X'\tht'_n G)(X'\sa) = X'G\tht_n$ in $\ulk{\Cl^+(\b\De_n^\#)}$.
\end{itemize}
\end{Remark}

{\it Proof.} Ad (1). We have 
\[
\barcl
(X\tht_n F)(X\rh)
& = & (X\tht_n F)(X\et^-\TTT F)(XF\sa^- F)(XF\TTT'\eps^-)        \\
& = & (X\tht_n F)([(X\et^-)^{+1}] F)(XF\sa^- F)([XF^{+1}]\eps^-) \\
& = & ([X\et^-]^{+1} F)(XFG\tht_n F)(XF\sa^- F)([XF^{+1}]\eps^-) \\
& = & ([X\et^- F]^{+1})(XFG\tht_n F)(XF\sa^- F)([XF^{+1}]\eps^-) \\
& = & ([XF\eps]^{+1})(XFG\tht_n F)(XF\sa^- F)([XF^{+1}]\eps^-)   \\
& = & XF\tht'_n \; .                                             \\
\ea
\]
Ad (2). We have
\[
\barcl
(X'\tht'_n G)(X'\sa)
& = & ([X'\eps]^{+1}G)(X'G\tht_n FG)(X'\sa^- FG)([X'^{+1}]\eps^- G)(X'\sa) \\
& = & ([X'\eps]^{+1}G)(X'G\tht_n FG)(X'\sa^- FG)([X'^{+1}]G\et)(X'\sa)     \\
& = & ([X'\eps]^{+1}G)(X'G\tht_n FG)([X'G^{+1}]\et)(X'\sa^-)(X'\sa)        \\
& = & ([X'\eps G]^{+1})(X'G\tht_n FG)([X'G^{+1}]\et)                       \\
& = & ([X' G\et^-]^{+1})(X'G\tht_n FG)([X'G^{+1}]\et)                      \\
& = & ([X' G]^{+1}\et^-)(X'G\tht_n FG)([X'G^{+1}]\et)                      \\
& = & (X'G\tht_n)([X' G^{+1}]\et^-)([X'G^{+1}]\et)                         \\
& = & X'G\tht_n\; .                                                        \\
\ea
\]
\qed

\subsection{Some lemmata}
\label{SecSomeLem}

Let $\El$ be a Frobenius category, let $\Bl\tm\El$ be its full subcategory of bijective objects; cf.\ e.g.\ \bfcite{Ku05}{Def.\ A.5}. We use the notations and conventions of \bfcite{Ku05}{\S A.2.3},
in particular those of \bfcite{Ku05}{Ex.\ A.6.(2)}.

Let $n\ge 0$. Let $E\tm\b\De_n^\#$ be a convex full subposet, i.e.\ whenever $\xi,\, \ze\,\in\, E$ and $\lm\in\b\De_n^\#$ such that $\xi\le\lm\le\ze$, then $\lm\in E$; cf.~\bfcite{Ku05}{\S 2.2.2.1}.
For instance, $\De_n^\trud\tm\b\De_n^\#$ is such a convex full subposet; cf.\ Notation~\ref{NotTDS1_5}.

A {\it pure square} in $\El$ is a commutative quadrangle $(A,B,C,D)$ with pure short exact diagonal sequence \mb{$(A,B\ds C,D)$}; cf.~\bfcite{Ku05}{\S A.4}. 

Denote by $\El^\Box(E) \tm \El(E)$ the full subcategory determined by
\[
\Ob\El^\Box(E) \; :=\; \left\{ \;\; X\in \Ob\El(E)\;\;\; :\; 
\mb{
\begin{tabular}{rl}
1) & $X_{\al/\al}$ is in $\Ob\Bl$ for all $\al\in \b\De_n$ such that $\al/\al\in E$, and       \\
   & $X_{\al^{+1}/\al}$ is in $\Ob\Bl$ for all $\al\in \b\De_n$ such that $\al^{+1}/\al\in E$. \\
2) & For all $\de^{-1}\le\al\le\be\le\ga\le\de\le\al^{+1}$ in $\b\De_n$\\
   & such that $\ga/\al$, $\ga/\be$, $\de/\al$ and $\de/\be$ are in $E$, \\
   & the quadrangle \\
   & $\xymatrix{
      X_{\ga/\be}\ar[r]^x                       & X_{\de/\be} \\
      X_{\ga/\al}\ar[r]_x\ar[u]^x\ar@{}[ur]|{\Box} & X_{\de/\al}\ar[u]_x \\      
      }$ \\
   & is a pure square. \\
\end{tabular}
}\right\} \;\; .
\]
A particular case of this definition has been considered in \bfcite{Ku05}{\S 4.1}.

\subsubsection{Cleaning the diagonal}

\begin{Lemma}
\label{LemCl1}
Suppose given $X\in\Ob\El^\Box(\b\De_n^\trud)$. Suppose given $\be\in\b\De_n$ such that $0\le\be\le 0^{+1}$. 

There exists $\w X\in\Ob\El^\Box(\b\De_n^\trud)$ such that the following conditions {\rm (1a,\,1b,\,2)} hold. 

\begin{itemize}
\item[{\rm (1a)}] We have $\w X_{\al/\al} = X_{\al/\al}$ for $0\le\al\le 0^{+1}$ such that $\al\ne\be$.
\item[{\rm (1b)}] We have $\w X_{\be/\be} = 0$.
\item[{\rm (2)}] There exists an isomorphism $\w X\lraiso X$ in $\ulEl(\b\De_n^\trud)$.
\end{itemize}
\end{Lemma}

{\it Proof.} Pars pro toto, we consider the case $n = 4$ and $\be = 2$. We display $X$ as follows.
\[
\xymatrix@R=7mm{
                &                                        &                                        &                                        &                                        & X_{0^{+1}/0^{+1}}   \\
                &                                        &                                        &                                        & X_{4/4}\ar[r]^x                        & X_{0^{+1}/4}\ar[u]^x\\
                &                                        &                                        & X_{3/3}\ar[r]^x                        & X_{4/3}\ar[r]^x\ar[u]^x\ar@{}[ur]|\Box & X_{0^{+1}/3}\ar[u]^x\\
                &                                        & X_{2/2}\ar[r]^x                        & X_{3/2}\ar[r]^x\ar[u]^x\ar@{}[ur]|\Box & X_{4/2}\ar[r]^x\ar[u]^x\ar@{}[ur]|\Box & X_{0^{+1}/2}\ar[u]^x\\
                & X_{1/1}\ar[r]^x                        & X_{2/1}\ar[r]^x\ar[u]^x\ar@{}[ur]|\Box & X_{3/1}\ar[r]^x\ar[u]^x\ar@{}[ur]|\Box & X_{4/1}\ar[r]^x\ar[u]^x\ar@{}[ur]|\Box & X_{0^{+1}/1}\ar[u]^x\\
X_{0/0}\ar[r]^x & X_{1/0}\ar[r]^x\ar[u]^x\ar@{}[ur]|\Box & X_{2/0}\ar[r]^x\ar[u]^x\ar@{}[ur]|\Box & X_{3/0}\ar[r]^x\ar[u]^x\ar@{}[ur]|\Box & X_{4/0}\ar[r]^x\ar[u]^x\ar@{}[ur]|\Box & X_{0^{+1}/0}\ar[u]^x\\
}
\]
Set $\w X$ to be the following diagram.
\[
\xymatrix@!@R=-8mm@C=4mm{
                &                                        &                                                           &                                                                                           &                                                                                                 & X_{0^{+1}/0^{+1}}                                   \\
                &                                        &                                                           &                                                                                           & X_{4/4}\ar[r]^x                                                                                 & X_{0^{+1}/4}\ar[u]^x                                \\
                &                                        &                                                           & X_{3/3}\ar[r]^x                                                                           & X_{4/3}\ar[r]^x\ar[u]^x\ar@{}[ur]|\Box                                                          & X_{0^{+1}/3}\ar[u]^x                                \\
                &                                        & 0\ar[r]                                                   & X_{3/2}\ar[r]^x\ar[u]^x\ar@{}[ur]|\Box                                                    & X_{4/2}\ar[r]^x\ar[u]^x\ar@{}[ur]|\Box                                                          & X_{0^{+1}/2}\ar[u]^x                                \\
                & X_{1/1}\ar[r]^x                        & X_{2/1}\ar[r]^(0.4){\smatez{x}{x}}\ar[u]\ar@{}[ur]|\Box   & X_{3/1}\dk X_{2/2}\ar[r]^{\smatzz{x}{0}{0}{1}}\ar[u]^{\rsmatze{x}{-x}}\ar@{}[ur]|\Box     & X_{4/1}\dk X_{2/2}\ar[r]^(0.45){\smatzz{x}{0}{0}{1}}\ar[u]^{\rsmatze{x}{-x}}\ar@{}[ur]|\Box     & X_{0^{+1}/1}\dk X_{2/2}\ar[u]^{\rsmatze{x}{-x}}     \\
X_{0/0}\ar[r]^x & X_{1/0}\ar[r]^x\ar[u]^x\ar@{}[ur]|\Box & X_{2/0}\ar[r]^(0.4){\smatez{x}{x}}\ar[u]^x\ar@{}[ur]|\Box & X_{3/0}\dk X_{2/2}\ar[r]^{\smatzz{x}{0}{0}{1}}\ar[u]^{\smatzz{x}{0}{0}{1}}\ar@{}[ur]|\Box & X_{4/0}\ds X_{2/2}\ar[r]^(0.45){\smatzz{x}{0}{0}{1}}\ar[u]^{\smatzz{x}{0}{0}{1}}\ar@{}[ur]|\Box & X_{0^{+1}/0}\ds X_{2/2}\ar[u]^{\smatzz{x}{0}{0}{1}} \\
}
\]
Using the Gabriel-Quillen-Laumon embedding theorem, we see that $\w X$ is actually an object of $\El^\Box(\b\De_n^\trud)$; cf.~\mb{\bfcite{Ku05}{\S A.2.2; Lem.\ A.11}}.

Since $X_{2/2}$ is bijective, inserting the zero morphism on all copies of $X_{2/2}$ and the identity on all other summands yields an 
isomorphism $\w X\lraiso X$ in $\ulEl(\b\De_n^\trud)$.\qed

\begin{Lemma}
\label{LemCl2}
Suppose given $X\in\Ob\El^\Box(\b\De_n^\trud)$.

There exists $X'\in\Ob\El^\Box(\b\De_n^\trud)$ such that the following conditions {\rm (1,\,2)} hold. 

\begin{itemize}
\item[{\rm (1)}] We have $X'_{\al/\al} = 0$ for all $\al\in\b\De_n$ such that $0\le\al\le 0^{+1}$. 
\item[{\rm (2)}] There exists an isomorphism $X'\lraiso X$ in $\ulEl(\b\De_n^\trud)$.
\end{itemize}
\end{Lemma}

{\it Proof.} This follows by application of Lemma~\ref{LemCl1} consecutively for $\be = 0$, $\be = 1$, \dots, $\be = 0^{+1}$.\qed

\subsubsection{Horseshoe lemma}
\label{SecHorseshoe}

Recall that $\Bl^{\text{ac}}$ denotes the category of purely acyclic complexes with entries in $\Bl$, i.e.\ of complexes with entries in $\Bl$ that decompose into pure short exact sequences in 
$\El$; cf.~\bfcite{Ku05}{\S A.2.3}.

Suppose given $Y\in\Ob\El$. An object $B$ of $\Bl^{\text{ac}}$ is called a {\it (both-sided) bijective resolution} of $Y$ if $Y$ is isomorphic to $\Img(B^0\lra B^1)$. Note that a bijective resolution
of a bijective object is split acyclic.

We have a full and dense functor (\footnote{A functor induced by $\h F$ will play the role of $F$ of Setup \ref{Setup2}; cf.\ \S\ref{SecConseq} below.})
\[
\ba{lcl}
\Bl^\text{ac} & \lraa{\h F} & \El                \\
B             & \lramaps    & \Img(B^0\lra B^1)  \\
\ea
\]
We make the additional convention that if the image factorisation of a pure morphism $d$ in $\El$ is chosen to be $d = \b d \dot d$, then we choose the image factorisation $-d = \b d (-\dot d)$ 
over the same image object.

Pointwise application yields a functor $\Bl^{\text{\rm ac}}(\b\De_n^\trud) \lraa{\h F} \El(\b\De_n^\trud)$, which is an abuse of notation.

Suppose given $X\in\Ob\El^\Box(\b\De_n^\trud)$ such that $X_{\al/\al} = 0$ for all $0\le\al\le 0^{+1}$. 

In particular, $X_{\be/\al}\lra X_{\ga/\al}\lra X_{\ga/\be}$ is a pure short exact sequence for \mb{$0\le\al\le\be\le\ga\le 0^{+1}$}.

Recall that for $n\in\b\De_n$, we have $n+1 = 0^{+1}$; cf.~\bfcite{Ku05}{\S 1.1}.

\begin{Lemma}
\label{LemHorseshoe}
Suppose given a bijective resolution $C_{\al+1/\al}$ of $X_{\al+1/\al}$ for all $\al\in\b\De_n$ such that $0\le\al\le n$. 

Then there exists $B\in\Ob(\Bl^{\text{\rm ac}})^\Box(\b\De_n^\trud)$ such that {\rm (1,\,2,\,3)} hold.

\begin{itemize}
\item[{\rm (1)}] We have $B\h F \iso X$ in $\El^\Box(\b\De_n^\#)$.
\item[{\rm (2)}] We have $B_{\al/\al} = 0$ for all $0\le\al\le 0^{+1}$.
\item[{\rm (3)}] We have $B_{\al+1/\al} = C_{\al+1/\al}$ for all $\al\in\b\De_n$ such that $0\le\al\le n$.
\end{itemize}
\end{Lemma}

\sbq
 If $n = 2$, and if we restrict to $\{ 1/0,\,2/0,\,2/1 \} \tm \b\De_2^\trud$, we recover the classical horseshoe lemma in its bothsided Frobenius category variant.
\seq

{\it Proof.} For $0\le\al\le n$, we denote $\left(C_{\al+1/\al}^0\lraepia{\b d} X_{\al+1/\al}\right) := \left(C_{\al+1/\al}^0\lraepi C_{\al+1/\al}\h F\lraiso X_{\al+1/\al}\right)$. 

By duality and by induction, it suffices to find a morphism $Y\lra X$ in $\El^\Box(\b\De_n^\trud)$ such that (i,\,ii,\,iii) hold.

\begin{itemize}
\item[(i)] We have $Y_{\be/\al}\in\Ob\Bl$ for all $0\le\al\le\be\le 0^{+1}$.
\item[(ii)] We have $Y_{\al/\al} = 0$ for all $0\le\al\le 0^{+1}$.
\item[(iii)] We have $(Y_{\al+1/\al}\lra X_{\al+1/\al}) \= (C_{\al+1/\al}^0\lraepia{\b d} X_{\al+1/\al})$ for all $\al\in\b\De_n$ such that $0\le\al\le n$.
\end{itemize}

Note that any morphism $Y\lra X$ fulfilling (i,\,ii,\,iii) consists pointwise of pure epimorphisms, and that the kernel of such a morphism $Y\lra X$ taken in $\El(\b\De_n^\trud)$ is in 
$\Ob\El^\Box(\b\De_n^\trud)$.

To construct $Y\lra X$, we let 
\[
Y_{\ga/\al}\; := \; \Ds_{\be\,\in\,\b\De_n,\; \al\le\be < \ga } C_{\be+1/\be}^0
\]
for $0\le\al\le\ga\le 0^{+1}$. For $\ga/\al \le \ga'/\al'$, the diagram morphism $Y_{\ga/\al}\lra Y_{\ga'/\al'}$ is stipulated to be identical on the summands $C_{\be+1/\be}^0$
with $\al'\le\be < \ga$ and zero elsewhere. This yields $Y\in\Ob\El^\Box(\b\De_n^\trud)$. 

Given $0\le\al\le\ga\le 0^{+1}$, we let $Y_{\ga/\al}\lra X_{\ga/\al}$ be defined as follows. For $0\le\be\le n$, we choose $Y_{\be+1/\be} \lraa{e} X_{\be+1/0}$ such that 
\[
(Y_{\be+1/\be} \lraa{e} X_{\be+1/0}\lraa{x} X_{\be+1/\be}) \= (Y_{\be+1/\be} \lraa{\b d} X_{\be+1/\be})\; .
\]
The component of the morphism
\[
(Y_{\ga/\al}\lra X_{\ga/\al}) \;\; =\;\; \left(\Ds_{\be\,\in\,\b\De_n,\; \al\le\be < \ga } C_{\be+1/\be}^0 \;\lra\; X_{\ga/\al}\right)
\]
at $\be$ is defined to be the composite
\[
(C_{\be+1/\be}^0 \lra X_{\ga/\al}) \;\; := \;\; (C_{\be+1/\be}^0 \lraa{e}  X_{\be+1/0} \lraa{x} X_{\ga/\al})\; .
\]
\qed

\subsubsection{Applying $\h F$ to a standard pure short exact sequence}

Recall that for $X\in\Ob\Bl^\text{\rm ac}$, we have chosen, in a functorial manner, a pure short exact sequence $X\lramono X\III\lraepi X\TTT$
with a bijective middle term $X\III$, where the letter $\III$ stands for ``injective''; cf.\ \bfcite{Ku05}{\S A.2.3}.

\begin{Lemma}
\label{LemFtoStand}
Suppose given $X\in\Ob\Bl^\text{\rm ac}$. There exists an isomorphism of pure short exact sequences in $\El$ as follows.
\[
\xymatrix@C=15mm{
X\h F\arm[r]\ar@{=}[d] & X\III\h F\are[r]\ar[d]_\wr & X\TTT\h F\ar@{=}[d] \\
X\h F\arm[r]           & X^1\are[r]                 & X\TTT\h F           \\
}
\]
Therein, the upper sequence results from an application of $\h F$ to the pure short exact sequence $X\lramono X\III\lraepi X\TTT$ in $\Bl^\text{\rm ac}$. The lower sequence is taken from the purely acyclic
complex $X$.
\end{Lemma}

{\it Proof.} Consider the following part of the pure short exact sequence $X\lramono X\III\lraepi X\TTT$ in $\Bl^\text{\rm ac}$; cf.~\bfcite{Ku05}{Ex.\ A.6}.
\[
\xymatrix@!@R=-2mm{
                                                                                                     & X^2\armfl{0.4}[rr]^(0.4){\smatez{1}{d}}                                                                 &                                                                                           & X^2\ds X^3\arefl{0.6}[rr]^(0.6){\rsmatze{-d}{1}}                                                                                                                         &                                                                     & X^3                                              \\
                                                                                                     &                                                                                                         &                                                                                           &                                                                                                                                                                          &                                                                     &                                                  \\
X\TTT\h F\arm[uur]^{\dot d}                                                                          &                                                                                                         &                                                                                           &                                                                                                                                                                          &                                                                     &                                                  \\
                                                                                                     &                                                                                                         &                                                                                           &                                                                                                                                                                          &                                                                     &                                                  \\
                                                                                                     & X^1\armfl{0.4}[rr]^(0.4){\smatez{1}{d}}\ar[uuuu]_d\are[uul]^*+<1mm,1mm>{\scm\b d}                       &                                                                                           & X^1\ds X^2\arefl{0.6}[rr]^(0.6){\rsmatze{-d}{1}}\ar[uuuu]_{\smatzz{0}{0}{1}{0}}                                                                                          &                                                                     & X^2\ar[uuuu]_{-d}                                \\
                                                                                                     &                                                                                                         &                                                                                           &                                                                                                                                                                          & X\III\h F\arm[ul]^{\dot\de}\are[d]^(0.3)*+<-1mm,0mm>{\scm\ka^-\b d} &                                                  \\
X\h F\arm[uur]^{\dot d}\armfl{0.4}[rr]_(0.4){\dot d}\arm[rrrru]^{\dot d \ka}|(0.57)\hole|(0.75)\hole &                                                                                                         & X^1\arm[uur]^{\smatez{1}{0}}\arefl{0.6}[rr]^(0.6){\b d\ru{-1}}\ar[urr]^(0.4)\ka_(0.4)\sim &                                                                                                                                                                          & X\TTT\h F\arm[uur]^{-\dot d}                                        &                                                  \\
                                                                                                     &                                                                                                         &                                                                                           &                                                                                                                                                                          &                                                                     &                                                  \\
                                                                                                     & X^0\armfl{0.4}[rr]^(0.4){\smatez{1}{d}}\ar'[uu][uuuu]_(0.4)d|(0.13)\hole\are[uul]^*+<1mm,1mm>{\scm\b d} &                                                                                           & X^0\ds X^1\arefl{0.6}[rr]^(0.6){\rsmatze{-d}{1}}\ar'[uu]_(0.7)*+<-1mm,0mm>{\smatzz{0}{0}{1}{0}}[uuuu]|(0.25)\hole\are[uul]^{\smatze{0}{1}}\are[uuur]_{\b\de}|(0.67)\hole &                                                                     & X^1\ar[uuuu]_{-d}\are[uul]^*+<1mm,1mm>{\scm\b d} \\
}
\]
We have added the image factorisations $X^0\lraepia{\b d} X\h F\lramonoa{\dot d} X^1$ and $X^0\dk X^1\lraepia{\b\de} X\III\h F\lramonoa{\dot\de} X^1\dk X^2$ 
of the respective differentials, resulting from an application of $\h F$. Factoring the differential of $X\TTT$ as
\[
\left(X^1\lraa{-d} X^2\right) \= \left(X^1\lraepia{\b d} X\TTT\h F\lramonoa{-\dot d} X^2\right)
\]
follows the additional convention made above.

Moreover, we have added the image factorisation $X^0\ds X^1\lraepifl{40}{\smatze{0}{1}} X^1\lramonofl{25}{\smatez{1}{0}} X^1\ds X^2$ and, accordingly, the isomorphism $X^1\lraisoa{\ka} X\III\h F$ that
satisfies $\ka\dot\de = \smatez{1}{0}$ and $\smatze{0}{1}\ka = \b\de$.

The horizontal pure short exact sequence $X\h F\lramonoa{\dot d} X^1\lraepia{\b d} X\TTT\h F$ lets all four arising parallelograms commute.

Now the sequence $X\lramono X\III\lraepi X\TTT$ maps to $X\h F\lrafl{25}{\dot d\ka} X\III\h F\lrafl{25}{\ka^-\b d} X\TTT\h F$, for the commutativities \mb{$(\dot d\ka)\dot\de = \dot d\smatez{1}{d}$} and 
$\b\de (\ka^- \b d) = \rsmatze{-d}{1}\b d$ hold.

In particular, the sequence $(X\lramono X\III\lraepi X\TTT)\h F$ actually is purely short exact.
\qed

\subsection{Stable vs.\ classically stable}

Let $\El$ be a Frobenius category, let $\Bl\tm\El$ be its full subcategory of bijective objects.

\subsubsection{$n$-triangles in the stable category}
\label{SecStableTriangles}

Recall that $\Bl^{\text{ac}}$ denotes the category of purely acyclic complexes with entries in $\Bl$; cf.~\S\ref{SecHorseshoe}.
Let $\Bl^{\text{sp\,ac}}\tm \Bl^{\text{ac}}$ denote the subcategory of split acyclic complexes.
Let $\uulEl = \Bl^{\text{ac}}/\Bl^{\text{sp\,ac}}$ denote the stable category of $\El$; cf.~\mb{\bfcite{Ku05}{Def.\ A.7}}. Let $\TTT$ be the automorphism on 
$\uulEl$ that shifts a complex to the left by one position, inserting signs; cf.~\bfcite{Ku05}{Ex.~A.6.(1)}. Then $(\uulEl,\TTT)$ carries a Heller triangulation $\tht$; 
cf.~\bfcite{Ku05}{Cor.\ 4.7}. In fact, we may, and will, choose the tuple of isotransformations $\tht = (\tht_n)_{n\ge 0}$ constructed in the proof of \bfcite{Ku05}{Th.\ 4.6}.

Suppose given $n\ge 0$ and $X\in\Ob(\Bl^{\text{ac}})^\Box(\b\De_n^\#)$; cf.~\bfcite{Ku05}{\S4.1} or \S\ref{SecSomeLem}. Now $X$ maps to an object $X\in\Ob\ulk{\uulEl^+(\b\De_n^\#)}$; 
cf.~\mb{\bfcite{Ku05}{Lem.\ A.29}}. Thus we have an isomorphism $[X]^{+1}\lrafl{25}{X\tht_n} [X^{+1}]$ in $\ulk{\uulEl^+(\b\De_n^\#)}$. By the construction in the proof of \bfcite{Ku05}{Th.\ 4.6}, 
there is a representative $[X]^{+1}\lraa{X\theta} [X^{+1}]$ in $\uulEl^+(\b\De_n^\#)$ of $X\tht_n$ such that in particular, there exists a morphism of pure short exact sequences
\[
\xymatrix@C=15mm{
X_{i/0}\armfl{0.35}[r]^(0.35){\smatez{x}{x}}\ar@{=}[d] & X_{i/i}\dk X_{0^{+1}\!/0}\arefl{0.6}[r]^(0.6){\rsmatze{x}{-x}}\ar[d] & X_{0^{+1}/i}\ar[d]^{X\h\theta_{i/0}} \\
X_{i/0}\armfl{0.35}[r]                                 & X_{i/0}\III\arefl{0.6}[r]                                            & X_{i/0}^{+1}                         \\
}
\]
for each $i\in [1,n]$; where $X\h\theta_{i/0}$ is a representative in $\Bl^{\text{ac}}$ for the morphism $X\theta_{i/0}$ in $\uulEl$; where the upper pure short exact sequence stems from the diagram $X$; 
and where the lower pure short exact sequence is the standard one as in \mb{\bfcite{Ku05}{Ex.\ A.6.(1)}}. In particular, $X\theta_{i/0}$ is an isomorphism in $\uulEl$.

Let $X^\tht\in\Ob\uulEl^{+,\,\per}(\b\De_n^\#)$ be defined as periodic prolongation of the image of the diagram $X|_{\b\De_n^\trud}$ in $\Ob\uulEl^+(\b\De_n^\trud)$ with $X_{0^{+1}/i}$ isomorphically 
replaced via $X\theta_{i/0}$ by $X_{i/0}^{+1}$ for all $i\in [1,n]$. For short, the rightmost column of the image of $X|_{\b\De_n^\trud}$ becomes standardised; cf.\ \bfcite{Ku05}{\S 2.1.3}. Using 
\[
\left.\left([X]^{+k}\lraiso [X^{+1}]^{+(k-1)}\lraiso\cdots\lraiso [X^{+k}]\right)\right|_{\b\De_n^\trud}
\]
for $k\ge 0$, and similarly for $k\le 0$, we obtain an isomorphism $X\smash{\lraa{\og}} X^\tht$ in $\uulEl^+(\b\De_n^\#)$ 
such that $\og_{i/0} = 1_{X_{i/0}}$ and $\og_{0^{+1}/i} = X\theta_{i/0}$ for $i\in [1,n]$; cf.\ Notation~\ref{NotTDS1_5}.(1).

\begin{Lemma}
\label{LemTDS5}
Given $n\ge 0$ and $X\in\Ob(\Bl^{\text{\rm ac}})^\Box(\b\De_n^\#)$, the periodic $n$\nbd-pretriangle $X^\tht$ is an $n$\nbd-triangle.
\end{Lemma}

\sbq
 The following proof is similar to the proof of Lemma~\ref{LemTDS2}.
\seq

{\it Proof.} We have to show that $X^\tht\tht_n = 1$; cf.~\bfcite{Ku05}{Def.\ 1.5.(ii.2)}. Since $\tht_n$ is a transformation, we have a commutative quadrangle
\[
\xymatrix@C=15mm{
[X^\tht]^{+1}\ar[r]^{X^\tht\tht_n}_\sim              & [(X^\tht)^{+1}]                 \\
[X]^{+1}\ar[r]^{X\tht_n}_\sim\ar[u]_{[\og]^{+1}}^\wr & [X^{+1}]\ar[u]_{[\og^{+1}]}^\wr \\
}
\]
in $\ulk{\uulEl^+(\b\De_n^\#)}\ru{-2.5}$. So we have to show that $(X\tht_n) [\og^{+1}] = [\og]^{+1}$. By \bfcite{Ku05}{Prop.\ 2.6}, it suffices to show that 
$(X\tht_n)|_{\dDe_n} [\og^{+1}]|_{\dDe_n} = [\og]^{+1}|_{\dDe_n\!}$. Now $[\og^{+1}]|_{\dDe_n} = 1_{[X^{+1}]|_{\dDe_n}}$ and, by construction, $(X\tht_n)|_{\dDe_n} = [\og]^{+1}|_{\dDe_n}$.\qed

\begin{Corollary}
\label{CorTDS5_5}
The Heller triangulated category $(\uulEl,\TTT,\tht)$ is closed.
\end{Corollary}

Cf.\ Definition~\ref{DefClosed}.

{\it Proof.} We can extend any morphism $X_{1/0}\lra X_{2/0}$ of $\Bl^\text{\rm ac}$ to an object of $(\Bl^{\text{\rm ac}})^\Box(\b\De_2^\#)$ by choosing $X_{1/0}\lramono X_{1/1}$ with $X_{1/1}$ bijective
and by choosing $X_{0^{+1}/0} = 0$, then forming pushouts, then choosing $X_{2/1}\lramono X_{2/2}$ with $X_{2/2}$ bijective, etc. Dually in the other direction. Then we apply Lemma~\ref{LemTDS5}.\qed

\subsubsection{The classical stable category under an additional hypothesis}

\paragraph{The hypothesis}\Absit
\label{SecHypo}

Let $\ulEl := \El/\Bl$ denote the classical stable category of $\El$.

Suppose given a set $\Dl$ of {\it distinguished} pure short exact sequences in $\El$ such that the following conditions hold.
\begin{itemize}
\item[(i)] The middle term of each distinguished pure short exact sequence is bijective.
\item[(ii)] For all $X\in\Ob\El$, there exists a unique distinguished pure short exact sequence with kernel term $X$.
\item[(iii)] For all $X\in\Ob\El$, there exists a unique distinguished pure short exact sequence with cokernel term $X$.
\end{itemize}

\paragraph{Consequences}\Absit
\label{SecConseq}

We shall define an endofunctor $\TTT'$ of $\ulEl$. 

{\it On objects.} Given $X\in\Ob\ulEl = \Ob\El$, there exists a unique distinguished pure short exact sequence with kernel term $X$. Let $X\TTT'$ be the cokernel term of this sequence. 

{\it On morphisms.} The image under $\TTT'$ of the residue class in $\ulEl$ of a morphism $X\lraa{f} Y$ in $\El$ is represented by the morphism 
$X\TTT'\lraa{g} Y\TTT'$ in $\El$ if there exists a morphism of distinguished pure short exact sequences as follows.
\[
\xymatrix{
X\armfl{0.48}[r]\ar[d]_f & B\arefl{0.42}[r]\ar[d] & X\TTT'\ar[d]^{g} \\
Y\armfl{0.48}[r]         & C\arefl{0.42}[r]       & Y\TTT'            \\
}
\]
Then $\TTT'$ is an automorphism of $\ulEl$; i.e.\ there exists an inverse $\TTT'^-$, constructed dually, such that $\TTT'\TTT'^- = 1_{\ulEl}$ and $\TTT'^-\TTT' = 1_{\ulEl}$.

As usual, we shall write $X^{+1} := X\TTT'$ for $X\in\Ob\ulEl$; etc.

The functor
\[
\barcl
\uulEl & \lraa{F} & \ulEl             \\
B      & \lramaps & \Img(B^0\lra B^1) \\
\ea
\]
induced by $\h F$ is an equivalence; cf.~\S\ref{SecHorseshoe}, \bfcite{Ku05}{Lem.\ A.1}. Splicing purely acyclic complexes from distinguished pure short exact sequences, we obtain an inverse equivalence 
$\uulEl\llaa{G}\ulEl$. Define $\TTT'G\lraa{\sa} G\TTT$ at $Y\in\Ob\ulEl = \Ob\El$ by letting
\[
(Y\sa^i) \; :=\; (-1_{(YG)^{i+1}})^i\; .
\]
Note that $Y\sa F = 1_{Y\TTT'GF} = 1_{YG\TTT F}$; cf.\ \S\ref{SecHorseshoe}.

Suppose given $Y\in\Ob\El$. We have a commutative diagram
\[
\xymatrix@C=4mm{
                      & (YG)^1                 &             \\ 
Y\ar[rr]_\sim\arm[ur] &                        & YGF\arm[ul] \\
                      & (YG)^0\are[ul]\are[ur] & ,           \\
}
\]
consisting of two image factorisations of the differential $(YG)^0 \lra (YG)^1$ and the induced isomorphism between the images $Y\lraiso YGF$ that makes the upper and the lower triangle 
in this diagram commute.

The residue class in $\ulEl$ of this induced morphism $Y\lraiso YGF$ shall be denoted by $Y\lraisoa{Y\eps} YGF$. Letting $Y$ vary, this gives rise to an isotransformation 
$1_{\ulEl}\lraisoa{\eps} GF$. Since $\eps$ is a transformation, we have $(Y\eps)(YGF\eps) = (Y\eps)(Y\eps GF)$, whence $GF\eps = \eps GF$. Thus there is an isotransformation $FG\lraisoa{\et} 1_{\uulEl}$ 
such that $(\eps G)(G\et) = 1_G$ and $(F\eps)(\et F) = 1_F$. Namely, for $\et$ we may take the inverse image under $F$ of $F\eps^-$.

So we are in the situation of Setup~\ref{Setup2} of \S\ref{SecDetecting}. Define $(\TTT F\lraa{\rh} F\TTT')$ as in \S\ref{SecDetecting}.

\begin{Proposition}
\label{PropTDS}\Absit
\begin{itemize}
\item[{\rm (1)}] By transport from $(\uulEl,\TTT,\tht)$ via $F$ and $G$, we obtain a closed Heller triangulation $\tht'$ on $(\ulEl,\TTT')$.
\item[{\rm (2)}] Suppose given $X'\in\Ob\,\ulEl^{+,\,\per}(\b\De_n^\#)$. Then $X'$ is an $n$\nbd-triangle if and only if $X'G^\sa$ is an $n$\nbd-triangle.
\item[{\rm (3)}] Suppose given $X\in\Ob\,\uulEl^{+,\,\per}(\b\De_n^\#)$. Then $X$ is an $n$\nbd-triangle if and only if $XF^\rh$ is an $n$\nbd-triangle.
\end{itemize}
\end{Proposition}

{\it Proof.} Assertion (1) follows by Lemmata~\ref{LemTDS1} and~\ref{LemTDS3}; cf.\ Corollary~\ref{CorTDS5_5}. Assertions (2,\,3) follow by Lemma~\ref{LemTDS2}.
\qed

\vs

Recall that $XF = X\h F$ in $\Ob\ulEl = \Ob\El$ for $X\in \Ob\uulEl = \Ob\Bl^\text{\rm ac}$.

\begin{Lemma}
\label{LemCalcRho}
Suppose given $X\in\Ob\Bl^\text{\rm ac}$. We have a morphism of pure short exact sequences
\[
\xymatrix@C=10mm{
XF\armfl{0.45}[r]\ar@{=}[d] & X^1\are[r]\ar[d] & X\TTT F\ar[d]       \\
XF\armfl{0.45}[r]           & (XFG)^1\are[r]   & XF\TTT'             \\
}
\]
in $\El$ such that its morphism $X\TTT F\lra XF\TTT'$ represents $X\rh$ in $\ulEl$. Here, the upper pure short exact sequence is taken from the purely acyclic complex $X$; the lower pure short exact 
sequence is distinguished.
\end{Lemma}

{\it Proof.} Given $X\in\Ob\Bl^\text{\rm ac}$, we can form a commutative diagram in $\El$ as follows.
\[
\xymatrix@!@C=-18mm@R=-22mm{
                       & X^2\ar[rr]                    &                                           & (XFG)^2                            &                               \\
X\TTT F\arm[ur]\ar[rr] &                               & XF\TTT'\arm[ur]\ar[rr]^(0.3){XF\TTT'\eps} &                                    & XF\TTT'GF = XFG\TTT F\arm[ul] \\
                       & X^1\ar[rr]\ar'[u][uu]\are[ul] &                                           & (XFG)^1\ar'[u][uu]\are[ul]\are[ur] &                               \\
XF\arm[ur]\ar@{=}[rr]  &                               & XF\arm[ur]\ar[rr]^(0.3){XF\eps}           &                                    & XFGF\arm[ul]                  \\
                       & X^0\ar[rr]\ar'[u][uu]\are[ul] &                                           & (XFG)^0\ar'[u][uu]\are[ul]\are[ur] &                               \\
}
\]
The morphisms $(XFG)^0\lraepi XF$, $XF\lramono (XFG)^1$, $(XFG)^1\lraepi XF\TTT'$ and $XF\TTT'\lramono (XFG)^2$ appear in distinguished pure short exact sequences.
Moreover, by abuse of notation, we have written $XF\eps$ resp.\ $XF\TTT'\eps$ for representatives in $\El$ of the respective morphisms in $\ulEl$. 

The partially displayed morphism of complexes $X\lra XFG$ represents $X\lrafl{25}{X\et^-} XFG$ in $\uulEl$, for $F$ maps the morphism represented by $X\lra XFG$ 
to $XF\eps = X\et^- F$.

Therefore, the composite morphism
\[
(X\TTT F\;\lra\; XF\TTT'\;\mrafl{25}{XF\TTT'\eps}\; XF\TTT' GF \= XFG\TTT F)
\]
from this diagram represents 
\[
(X\;\lrafl{25}{X\et^-}\; XFG)\TTT F\; ;
\]
note that there are no signs to be inserted at the respective pure epimorphisms of the image factorisation chosen by $F$; cf.\ \S\ref{SecHorseshoe}.
Thus the morphism $X\TTT F\lra X F\TTT'$ from this diagram represents 
\[
\big(X\TTT F\;\mrafl{25}{X\et^-\TTT F}\; XFG\TTT F \= XF\TTT' GF \;\mrafl{25}{XF\TTT'\eps^-}\; XF\TTT'\big) \;\;\=\;\; \big(X\TTT F\;\lrafl{25}{X\rh} XF\TTT'\big)\; .
\]
\qed

\paragraph{Standardisation by substitution of the rightmost column}\Absit
\label{SecStandardise}

\sbq
 We mimic the construction $X\ramaps X^\tht$ made in~\S\ref{SecStableTriangles}, now for $\El$ instead of $\Bl^\text{ac}$.
\seq

Denote by $\El^\Box(\b\De_n^\trud) \lraa{M'} \ulEl^+(\b\De_n^\trud)$ the residue class functor; cf.~\S\ref{SecSomeLem}, \bfcite{Ku05}{\S 2.1.3, \S 4.1}. 
Denote by $\ulEl^{+,\,(\Box)}(\b\De_n^\trud)$ the full subcategory of $\ulEl^+(\b\De_n^\trud)$ whose set of objects is given by
\[
\Ob\ulEl^{+,\,(\Box)}(\b\De_n^\trud) \; := \; \big(\Ob\El^\Box(\b\De_n^\trud)\big)M'\; .
\]
So $\ulEl^{+,\,(\Box)}(\b\De_n^\trud)$ is defined to be the ``full image'' in $\ulEl^+(\b\De_n^\trud)$ of the residue class functor $M'$.

Suppose given $n\ge 0$ and $X\in\Ob\El^\Box(\b\De_n^\trud)$. Recall that $\dDe_n = [1,n]$ is identified with $\{i/0 \; :\; i\in [1,n]\} \tm \b\De_n^\trud$.
Write $X_{\ast/0} := X|_{\dDe_n} = (i\lramaps X_{i/0}) \in\Ob\El(\dDe_n)$ and $X_{0^{+1}/\ast} := (i\lramaps X_{0^{+1}/i}) \in\Ob\El(\dDe_n)$, analogously for morphisms; analogously for
objects in $\ulEl^+(\b\De_n^\trud)$ and their morphisms.

Let the isomorphism $\smash{X_{0^{+1}/\ast}\lraisofl{22}{X\tau} X_{\ast/0}^{+1}}$ in $\ulEl(\dDe_n)$ be defined by morphisms of pure short exact sequences
\[
\xymatrix@C=15mm{
X_{i/0}\armfl{0.35}[r]^(0.35){\smatez{x}{x}}\ar@{=}[d] & X_{i/i}\dk X_{0^{+1}\!/0}\arefl{0.6}[r]^(0.6){\rsmatze{x}{-x}}\ar[d] & X_{0^{+1}/i}\ar[d]^{X\h\tau_i} \\
X_{i/0}\armfl{0.35}[r]                                 & B_{i/0}\arefl{0.6}[r]                                                & X_{i/0}^{+1}                   \\
}
\]
for $i\in [1,n]$; where $X\h\tau_i$ is a representative in $\El$ for the morphism $X\tau_i$ in $\ulEl$; where the upper pure short exact sequence stems from the diagram $X$; 
and where the lower pure short exact sequence is distinguished.

In this way, we get an isotransformation $\tau$ between the functors $(-)_{0^{+1}/\ast}$ and $(-)_{\ast/0}^{+1}$ from $\ulEl^{+,\,(\Box)}(\b\De_n^\trud)$ to $\ulEl(\dDe_n)$.

Let $\ulEl^{+,\,\per}(\b\De_n^\trud)$ be the (in general not full) subcategory of $\ulEl^+(\b\De_n^\trud)$ given by the set of objects
\[
\Ob\ulEl^{+,\,\per}(\b\De_n^\trud) \; :=\; \left\{ Y\in\Ob\ulEl^+(\b\De_n^\trud) \; :\; \text{$Y_{0^{+1}/\ast} = Y_{\ast/0}^{+1}$ in $\Ob\ulEl(\dDe_n)$} \right\}\; ,
\]
and by the set of morphisms
\[
\liu{\ulEl^{+,\,\per}(\b\De_n^\trud)}{(Y,Y')} \= \{ f\in \liu{\ulEl^+(\b\De_n^\trud)}{(Y,Y')} \; :\; \text{$f_{0^{+1}/\ast} = f_{\ast/0}^{+1}$ in $\ulEl(\dDe_n)$}\}\; ,
\]
for $Y,\, Y'\,\in\,\Ob\ulEl^{+,\,\per}(\b\De_n^\trud)$.

Given $X\,\in\,\Ob\El^\Box(\b\De_n^\trud)$, we let $X^\tau\in\Ob\ulEl^{+,\,\per}(\b\De_n^\trud)$ be defined as the diagram $X$ with $X_{0^{+1}/\ast}$ isomorphically replaced via $X\tau$ by 
$X_{\ast/0}^{+1}$. For short, the rightmost column of $X$ becomes standardised to obtain $X^\tau$.

Given $X,\, X'\,\in\,\Ob\El^\Box(\b\De_n^\trud)$, a morphism $X\lraa{f} X'$ in $\ulEl^+(\b\De_n^\trud)$ induces a morphism $X^\tau \lraa{f^\tau} X'^\tau$ in 
$\ulEl^{+,\,\per}(\b\De_n^\trud)$. Namely, we let $f^\tau_{\be/\al} := f_{\be/\al}$ for $0\le\al\le\be\le n$, and we let $f^\tau_{0^{+1}/\ast}$ be characterised by the commutative quadrangle
\[
\xymatrix{
X_{0^{+1}/\ast}\ar[r]^{X\tau}_\sim\ar[d]_{f_{0^{+1}/\ast}}  & X_{\ast/0}^{+1} \ar[d]^{f^\tau_{0^{+1}/\ast}} \\
X'_{0^{+1}/\ast}\ar[r]^{X'\tau}_\sim                        & {X'}_{\ast/0}^{+1}                            \\
}
\]
in $\ulEl(\dDe_n)$. In particular, since $\tau$ is an isotransformation, we have $f^\tau_{0^{+1}/\ast} = (f^\tau_{\ast/0})^{+1}$.

\begin{Remark}
\label{RemTDS6}
The constructions made above define a functor 
\[
\barcl
\ulEl^{+,\,(\Box)}(\b\De_n^\trud) & \lrafl{25}{(-)^\tau} & \ulEl^{+,\,\per}(\b\De_n^\trud) \\
                               X  & \lramapsfl{}{}       & X^\tau \; .                     \\
\ea
\]
\end{Remark}

\paragraph{$n$-triangles in the classical stable category}\Absit
\label{SecTrianglesClassical}

\begin{Proposition}
\label{PropTDS7}
Suppose given $n\ge 0$ and $X\in\Ob\El^\Box(\b\De_n^\trud)$. 

The periodic prolongation of \mb{$X^\tau\in\Ob\ulEl^{+,\,\per}(\b\De_n^\trud)$} to an object of $\ulEl^{+,\,\per}(\b\De_n^\#)$ is an $n$\nbd-triangle with respect to the triangulation $\tht'$ on 
$(\ulEl,\TTT')$ obtained as in Proposition {\rm\ref{PropTDS}}.
\end{Proposition}

{\it Proof.} By Lemma~\ref{LemCl2}, there exists $X'\in\Ob\El^\Box(\b\De_n^\trud)$ such that $X'_{\al/\al} = 0$ for all $0\le\al\le 0^{+1}$ and such that $X$ is isomorphic to $X'$ in 
$\ulEl^+(\b\De_n^\trud)$. By Remark~\ref{RemTDS6}, the object $X^\tau$ is isomorphic to $X'^\tau$ in $\ulEl^{+,\,\per}(\b\De_n^\trud)$. Thus the periodic prolongation of 
$X^\tau$ is an $n$\nbd-triangle if and only if that of $X'^\tau$ is; cf.~\bfcite{Ku05}{Lem.\ 3.4.(4)}.

Therefore, we may assume that $X_{\al/\al}\iso 0$ for all $0\le\al\le 0^{+1}$.

Let $\w X\in\Ob(\Bl^{\text{\rm ac}})^\Box(\b\De_n^\trud)$ be such that there exists an isomorphism $\smash{X \lraisoa{\h a} \w X\h F}$ in $\El^\Box(\b\De_n^\trud)$ and such that 
$\w X_{\al/\al} = 0$ for all $0\le\al\le 0^{+1}$; cf.\ Lemma~\ref{LemHorseshoe}.
Denote by $\smash{X \lraisoa{a} \w X F}\ru{4}$ the isomorphism in $\ulEl^+(\b\De_n^\trud)$ represented by $\smash{X \lraisoa{\h a} \w X\h F}\ru{4}$.

Let $\w{\w X}\in\Ob(\Bl^{\text{\rm ac}})^\Box(\b\De_n^\#)$ be such that $\w{\w X}|_{\b\De_n^\trud} = \w X$. 
By Lemma~\ref{LemTDS5}, the periodic $n$\nbd-pretriangle $\w{\w X}^\tht\in\Ob\uulEl^{+,\,\per}(\b\De_n^\#)$ is an $n$\nbd-triangle. Note that $\w{\w X}^\tht\ru{5}$ depends only on $\w X$, not on the 
choice of $\w{\w X}$.

Thus, by Proposition~\ref{PropTDS}.(3), $\w{\w X}^\tht F^\rh \in\Ob\ulEl^{+,\,\per}(\b\De_n^\#)$ is an $n$\nbd-triangle. Therefore, it suffices to show that $X^\tau$ and 
$\w{\w X}^\tht F^\rh|_{\b\De_n^\trud}$ are isomorphic in $\ulEl^{+,\,\per}(\b\De_n^\trud)$, for then their periodic prolongations are isomorphic in $\ulEl^{+,\,\per}(\b\De_n^\#)$, which in turn shows 
the periodic prolongation of $X^\tau$ to be an $n$\nbd-triangle; cf.~\bfcite{Ku05}{Lem.\ 3.4.(4)}.

We have a composite isomorphism
\[
X^\tau\;\llaiso\; X\;\lraisoa{a}\; \w X F \= \w{\w X} F|_{\b\De_n^\trud} \;\lraiso\; \w{\w X}^\tht F|_{\b\De_n^\trud} \;\lraiso\; \w{\w X}^\tht F^\rh|_{\b\De_n^\trud}
\]
in $\ulEl^+(\b\De_n^\trud)$. We {\it claim} that it lies in $\ulEl^{+,\,\per}(\b\De_n^\trud)$.

Suppose given $i\in [1,n]$. On $i/0$, this composite equals $a_{i/0}$. Thus we have to show that on $0^{+1}/i$, this composite equals $a_{i/0}^{+1} = a_{i/0}\TTT'$. Consider, to this end, the 
following morphisms of pure short exact sequences in $\El$.
\[
\xymatrix@C=12mm@R=6mm{
X_{i/0}\armfl{0.45}[r]                                                                 & (X_{i/0}G)^1\are[r]                                                             & X_{i/0}\TTT'                                                             \\
X_{i/0}\armfl{0.45}[r]^(0.46){x\ru{-0.5}}\ar@{=}[u]\ar[d]^(0.45){\h a_{i/0}}_(0.45)\wr & X_{0^{+1}/0}\are[r]^{-x\ru{-0.5}}\ar[u]\ar[d]^(0.45){\h a_{0^{+1}/0}}_(0.45)\wr & X_{0^{+1}/i}\ar[u]_{X\h\tau_i}\ar[d]^(0.45){\h a_{0^{+1}/i}}_(0.45)\wr   \\
\w X_{i/0}\h F\armfl{0.45}[r]^(0.46){\w x\h F}\ar@{=}[d]                               & \w X_{0^{+1}/0}\h F\are[r]^{-\w x\h F}\ar[d]                                    & \w X_{0^{+1}/i}\h F\ar[d]^{X\h\theta_{i/0}\h F}                          \\
\w X_{i/0}\h F\armfl{0.45}[r]\ar@{=}[d]                                                & \w X_{i/0}\III\h F\are[r]\ar[d]                                                 & \w X_{i/0}\TTT\h F\ar@{=}[d]                                             \\
\w X_{i/0}\h F\armfl{0.45}[r]\ar@{=}[d]                                                & (\w X_{i/0})^1\are[r]\ar[d]                                                     & \w X_{i/0}\TTT\h F\ar[d]                                                 \\
\w X_{i/0}\h F\armfl{0.45}[r]                                                          & (\w X_{i/0}\h F G)^1\are[r]                                                     & \w X_{i/0}\h F\TTT'                                                      \\
}
\]
The fourth sequence is purely short exact by Lemma~\ref{LemFtoStand}.

The first morphism from above arises by definition of $\h\tau_i$; cf.~\S\ref{SecStandardise}. The second morphism is taken from $\h a$. The third morphism arises by definition of $\h\theta_{i/0}$ and 
an application of $\h F$; cf.~\S\ref{SecStableTriangles}. The fourth morphism is given by Lemma~\ref{LemFtoStand}. The fifth morphism is given by Lemma~\ref{LemCalcRho}.

The first and the sixth pure short exact sequence are distinguished, and so the {\it claim} and hence the proposition follow. 

\end{footnotesize}

\parskip0.0ex
\begin{footnotesize}

\parskip1.2ex

\vspace*{5mm}

\begin{flushright}
Matthias K\"unzer\\
Lehrstuhl D f\"ur Mathematik\\
RWTH Aachen\\
Templergraben 64\\
D-52062 Aachen \\
kuenzer@math.rwth-aachen.de \\
www.math.rwth-aachen.de/$\sim$kuenzer\\
\end{flushright}
\end{footnotesize}

\end{document}